\documentclass[11pt]{amsart}
\usepackage{amsmath,amssymb,amsthm,amscd,latexsym}
\usepackage[all]{xy}
\usepackage[dvips]{graphicx}
\usepackage{hyperref}

\input xypic
\xyoption{all}

\newcommand{\rar}{\rightarrow}
\newcommand{\lar}{\longrightarrow}
\newcommand{\llar}{-\kern-5pt-\kern-5pt\longrightarrow}

\newtheorem{Theorem}{Theorem}[section]
\newtheorem{Lemma}[Theorem]{Lemma}

\newtheorem{Proposition}[Theorem]{Proposition}
\newtheorem{Remark}[Theorem]{Remark}

\def\sqr#1#2{{\vcenter{\hrule height.#2pt
        \hbox{\vrule width.#2pt height#1pt \kern#1pt
            \vrule width.#2pt}
        \hrule height.#2pt}}}
\def\phi{\varphi}
\def\demo{\noindent{\bf Proof. }}
\def\square{\mathchoice\sqr64\sqr64\sqr{4}3\sqr{3}3}
\def\qed{\hspace*{\fill} $\square$}

\def\xx{{\bf x}}

\def\fm{{\mathfrak m}}

\def\ker{{\rm ker}\,}

\def\restr{{\kern-1pt\restriction\kern-1pt}}

\begin{document}
\title[On a conjecture of Vasconcelos]{On a conjecture of Vasconcelos}

\author{Ricardo Burity, Aron Simis}


\author{Stefan O. Toh\v aneanu}

 \keywords{Rees algebra, Sylvester forms, almost Cohen--Macaulay, reduction number, monomials, initial ideal. \\
\indent The second author is supported by a CNPq grant and a PVNS Fellowship from CAPES. He thanks the Departamento de Matem\'atica of the Universidade Federal da Paraiba for providing an appropriate environment for discussions on this work.\\
\indent Burity address: Departamento de Matem\'atica, Universidade Federal da Paraiba, 58051-900 Jo\~ao Pessoa, PB, Brazil, email: ricardoburity@hotmail.com\\
\indent Simis address: Departamento de Matem\'atica, Universidade Federal de Pernambuco, 50740-640 Recife, PE, Brazil, email: aron@dmat.ufpe.br\\
\indent Tohaneanu  address: Department of Mathematics, University of Idaho, Moscow, Idaho 83844-1103, USA, email: tohaneanu@uidaho.edu.} 





\begin{abstract}
One studies the structure of the Rees algebra of an almost complete intersection monomial ideal of finite co-length in a polynomial ring over a field, assuming that the least pure powers of the variables contained in the ideal have the same degree. It is shown that the Rees algebra has a natural quasi-homogeneous structure and its presentation ideal is generated by explicit Sylvester forms. A consequence of these results is a proof that the Rees algebra is almost Cohen--Macaulay, thus answering affirmatively an important case of a conjecture of W. Vasconcelos.

\end{abstract}

\maketitle


\section*{Introduction}

Let $R:=k[x_1,\ldots,x_n]$ denote a polynomial ring over a field $k$.
In 2013 W. Vasconcelos formulated the conjecture that the Rees algebra of an Artinian almost complete intersection $I\subset R$ generated by monomials is almost Cohen-Macaulay. For the binary case (i.e., for $n=2$) a result of M. Rossi and I. Swanson (\cite[Proposition 1.9]{RoSw}) gives an affirmative answer to the conjecture. Recently, different proofs were established in the binary case of monomials of the same degree  as a consequence of work by T. Benitez and C. D'Andrea (\cite{Dandrea}), and independently, of work by the present second and third authors (\cite{sylStefan}).

One main tool in the binary case of an Artinian almost complete intersections $I$ of forms of the same degree is birationality.
The other two tools are the Ratliff--Rush filtration theory and the Huckaba--Marley criterion using a minimal reduction of $I$. While  the Ratliff-Rush filtration  gives no insight into the conjectured property of the Rees algebra beyond the two variables case, using the criterion of Huckaba--Marley, would probably require as much calculation and besides lead one into no reasonable bound to manage the partial lengths.
We add the fact that even when the uniformity assumption degenerates into equigrading, birationality for more than two variables is an issue, and hence computing the first Hilbert coefficient of $R/I$ becomes a hardship.

\smallskip

The purpose of this paper is to tackle the case of arbitrary number of variables with an extra condition on the degrees of the monomials, called {\em uniformity}.
Under this condition, we answer affirmatively the stated conjecture. In our opinion this contributes a significant step toward the general case, since we have in mind a couple of procedures to reducing the case of general exponents to this one. Notwithstanding the seemingly simple case of a monomial ideal, as compared to the problem of ideals generated by arbitrary forms -- a situation still lacking a bona fide conjecture -- its general case may require ahead an additional tour de force beyond the facilitation provided by the methods of the present paper.

The method in the present paper emphasizes the structure of the presentation ideal of $\mathcal{R}_R(I)$ that may benefit from the appeal to Sylvester forms, as we understand them in their modern algebraic formulation.
However, additional muscle work soon became imperative, so the technology comes from three  sources: first, a thorough employ of the natural quasi-homogeneous grading over $k$ of the presentation ideal of  $\mathcal{R}_R(I)$, compatible with the usual standard grading of $\mathcal{R}_R(I)$ over $R$; second, mastering the overwhelming presence of a sequence of iterated Sylvester forms that are Rees generators; third,  the perception of an underlying monomial order in the ambient polynomial ring $S\supset R$ of the presentation of $\mathcal{R}_R(I)$, allowing for a careful computation of certain colon ideals crucial for extracting the homological nature of $\mathcal{R}_R(I)$.
To our knowledge, the systematic role of these three tools has not been sufficiently emphasized elsewhere.

The core of the paper is confined to the first two sections, while we collected the technically hard proofs in the third section in order to avoid distraction.
In the first section one develops the details of a very precise set of generators of the Rees presentation ideal, drawing upon a weighted grading naturally stemming from the form of the monomial generators of $I$.
Of course, it is well-known that ideals of relations of monomials are generated by binomials. However, for the sake of an efficient generation we need to show that these binomials acquire a special form due to the nature of the given monomials.

One shows that the relation type of $I$ equals the reduction number of $I$ plus $1$ and, moreover, state a precise count of the number of the generators in each external (i.e., presentation) degree.
Finally, one dedicates a stretch of the section to the identification of these binomial generators as iterated Sylvester forms.

In the subsequent section one states that the above generators can be ordered in a such a way as to describe the Rees presentation ideal $\mathcal{I}$ of $I$ by a finite series of subideals of which any two consecutive ones have a monomial colon ideal.
By inducting on the length of this series one is then able to consider mapping cones iteratively culminating with $\mathcal{I}$ itself. As a consequence,   the Rees algebra $\mathcal{R}_R(I)$ will be almost Cohen--Macaulay, thus answering affirmatively in this case a conjecture of Vasconcelos stated in \cite[Conjecture 4.15]{syl3}.

The preliminaries of this section require dealing at length with initial ideals and their colon ideals.
The calculations along this line of approach though basically straightforward are quite lengthy and seem to be unavoidable.
For the purpose of not disturbing the readership smoothness of the main results, we collected those proofs in the subsequent section.
Although the details of the proofs can be avoided in a first reading, they constitute the fine tissue legitimating the main results of the paper.

The main results of this paper are Theorem~\ref{format}, Theorem~\ref{gens_relation_type},
Proposition~\ref{colon-of-initials}, Proposition~\ref{successive_colon} and Theorem~\ref{rees_is_acm}.

\medskip

A substantial portion of the computational results are part of the PhD thesis of the first author.

\section{Efficient generation}

Let $R:=k[x_1,\ldots,x_n]$ denote a polynomial ring over a field $k$.
Given integers $0<b<a$, the monomial ideal   $I:=(x_1^{a},\ldots,x_n^{a}, (x_1\cdots x_n)^b)\subset R$ will be called {\em uniform}.

Our main focus is the presentation of the Rees algebra $\mathcal{R}_R(I)$ over a polynomial ring $S:=R[y_1,\ldots,y_n,w]$:
$$\mathcal I:=\ker(S \lar R[It]),\;y_j\mapsto x_j^at, w\mapsto (x_1\cdots x_n)^bt.$$
The presentation ideal  $\mathcal{I}\subset S$ is often referred to as the {\em Rees ideal} of $I$ and $y_1,\ldots,y_n,w$ as the {\em presentation} or {\em external} variables.
We will moreover let ${\mathcal L}\subset \mathcal{I}$ denote the set of generators coming from the syzygies of $I$.

A major question is a lower bound for the depth of $\mathcal{R}_R(I)$,  where the depth is computed on the maximal graded ideal $(\fm, S_+)$, with $\fm=(x_1,\ldots,x_n)$.
Knowingly, $\mathcal{R}_R(I)$ is Cohen--Macaulay when its depth attains the maximum value in the inequality
${\rm depth}(\mathcal{R}_R(I))\leq \dim\mathcal{R}_R(I)=n+1$.
One says that $\mathcal{R}_R(I)$ is {\em almost Cohen--Macaulay}  if ${\rm depth}(\mathcal{R}_R(I))\geq n$, a condition
 equivalent to $\mathcal{R}_R(I)$ having homological dimension $\leq n+1$ over $S$.

\medskip

As prolegomena, we restate the following valuable piece of information about the reduction number of the ideal $I$ proved in \cite{sylStefan}:

\begin{Proposition}\label{red_numbers} {\rm (\cite[Proposition 2.13]{sylStefan})}
For a uniform monomial ideal as above the following hold:
\begin{itemize}
\item[{\rm (a)}] $J:=(x_1^a,\ldots,x_n^a)$ is a minimal reduction of $I$ if and only if $nb\geq a${\rm ;} in this case, letting $1\leq p\leq n$ be the smallest integer such that $pb\geq a$ {\rm (}hence $(p-1)b<a${\rm )}, one has ${\rm red}_J(I)=p-1$.
\item[{\rm (b)}] If $nb<a$, then
$Q:=(x_1^a-x_n^a,\ldots,x_{n-1}^a-x_n^a,(x_1\cdots x_n)^b)$ is a minimal reduction of $I$ and ${\rm red}_Q(I)=n-1$.
\end{itemize}
\end{Proposition}

In particular, for $a\leq 2b$ the ideal  $J$ is a minimal reduction of $I$ with reduction number $1$, hence is $\mathcal{R}_R(I)$ is Cohen--Macaulay as is well-known.
Since this situation has no interest in our discussion, we will assume $a>2b$ throughout the paper.

\smallskip

In this part we search for a set of binomials of a particular form that minimally generate the Rees ideal $\mathcal{I}$ of $I$.
As we will contend in Theorem~\ref{format}, the ring $S$ admits a weighted grading under which an extra behavior will emerge.
For now, as a preamble we can prove a basic result that depends solely on the standard grading of $S$ as a polynomial ring over $R$.
It is well-known that ideals of relations of monomials are generated by binomials. In the present case, we show that these binomials acquire a special form due to the nature of the given monomials.
This step will be crucial in the subsequent unfolding.

\begin{Lemma}\label{basic_form}
Any binomial in $\mathcal I$ belonging to a set of minimal generators thereof is of the form
$$\mathbf{m}(\xx)w^{\delta}-\mathbf{n}(\xx)y_{i_1}^{\alpha_{i_1}}\cdots y_{i_s}^{\alpha_{i_s}},$$
 where  $\mathbf{m}(\xx), \mathbf{n}(\xx)$ are relatively prime monomials in $\xx=x_1,\ldots,x_n$ and  $1\leq i_1<\cdots < i_s\leq n$, $\alpha_{i_j}>0$.
\end{Lemma}
\demo One has to show that, for no $1\leq i\leq n$ do $y_i$ and $w$ divide the same monomial in the expression of a generating binomial.

Assuming the contrary, one has the following two possibilities for a binomial relation:

\smallskip

{\sc Case 1.} $y_1^{\alpha_1}\cdots y_t^{\alpha_t}w^{\delta}-x_1^{d_1}\cdots x_t^{d_t}x_{t+1}^{d_{t+1}}\cdots x_n^{d_n}y_{t+1}^{\alpha_{t+1}}\cdots y_n^{\alpha_n},$
where $\delta>0$ and  $\alpha_1,\ldots,\alpha_t\geq 1$.

\smallskip

Because of the homogeneity of the variables $y_1,\ldots,y_n,w$ and since  upon evaluation the degrees of $x_1,\ldots,x_n$ must match on the two sides, we obtain the numerical equalities

\begin{eqnarray}
\alpha_1+\cdots+\alpha_t+\delta&=&\alpha_{t+1}+\cdots+\alpha_n\nonumber\\
a\alpha_j+\delta b&=&d_j,j=1,\ldots,t\nonumber\\
\delta b&=&a\alpha_k+d_k,k=t+1,\ldots,n.\nonumber
\end{eqnarray}

From the first of these equalities we can assume that $\alpha_{t+1}\geq 1$ and, from the second one, that $d_1>a$. Then the binomial can be written as $$y_1^{\alpha_1}\cdots y_t^{\alpha_t}w^{\delta}- \underbrace{(K_{1,t+1}+x_{t+1}^ay_1)}_{x_1^ay_{t+1}}x_1^{d_1-a}\cdots x_t^{d_t}x_{t+1}^{d_{t+1}}\cdots x_n^{d_n}y_{t+1}^{\alpha_{t+1}-1}\cdots y_n^{\alpha_n},$$ where $K_{i,j}=x_{i}^{a}y_{j}-x_{j}^{a}y_{i},~i,j\in\{1,\ldots,n\},~i<j$.

Since $\mathcal I$ is a prime ideal, simplifying by $y_1$ due to  minimality, one obtains a binomial in $\mathcal I$ of the same shape with $y_1$ raised to the power $\alpha_1-1$. Iterating, we can replace the given generator by another one of the same shape, where the exponent of $y_1$ vanishes. But this contradicts the assumption that this exponent is nonzero.

\medskip

{\sc Case 2.} $x_1^{d_1}\cdots x_m^{d_m} y_{m+1}^{\alpha_{m+1}}\cdots y_t^{\alpha_t}w^{\delta}-x_{m+1}^{d_{m+1}}\cdots x_t^{d_t}y_1^{\alpha_1}\cdots y_m^{\alpha_m}x_{t+1}^{d_{t+1}}\cdots x_n^{d_n}y_{t+1}^{\alpha_{t+1}}\cdots y_n^{\alpha_n},$
where $\delta>0$ and  $\alpha_{m+1},\ldots,\alpha_t\geq 1$.

\smallskip

As before, one has the following equalities between the exponents:

\begin{eqnarray}
d_i+\delta b&=&\alpha_ia, i=1,\ldots,m\nonumber\\
a\alpha_j+\delta b&=&d_j,j=m+1,\ldots,t\nonumber\\
\delta b&=&a\alpha_k+d_k,k=t+1,\ldots,n.\nonumber
\end{eqnarray}

As $\delta>0$, the first set of equations gives $\alpha_1,\ldots,\alpha_m\geq 1$. The assumption $\alpha_j\geq 1, j=m+1,\ldots,t$, and the second set of equations give $d_{m+1},\ldots,d_t>a$. Then the binomial can be written in the form

{\small
$$x_1^{d_1}\cdots x_m^{d_m} y_{m+1}^{\alpha_{m+1}}\cdots y_t^{\alpha_t}w^{\delta}-\underbrace{(K_{m+1,1}+x_1^ay_{m+1})}_{x_{m+1}^ay_1}x_{m+1}^{d_{m+1}-a}\cdots x_t^{d_t}y_1^{\alpha_1-1}\cdots y_m^{\alpha_m}x_{t+1}^{d_{t+1}}\cdots x_n^{d_n}y_{t+1}^{\alpha_{t+1}}\cdots y_n^{\alpha_n}.$$
}

By the same token as above, one obtains a binomial in $\mathcal I$ of the same shape with  $y_{m+1}$ raised  to $\alpha_{m+1}-1$. Iterating on $\alpha_{m+1}$ as in the first case gives a contradiction -- note that, because $\alpha_1$ also drops  by $1$,  the first case is around the corner in the inductive process.
\qed

\bigskip

The following notation will be used throughout the rest of the paper: if $\{i_1,\ldots ,i_j\}$ is a subset of $\{1,\ldots,n\}$  we denote by $P(i_1,\ldots,i_j)$ the product of the variables belonging to the subset $\{x_1,\ldots,x_n\}\setminus \{x_{i_1},\ldots, x_{i_j}\}$.
A few times around we may deal with a similar situation where we may wish to stress that $\{i_1,\ldots ,i_j\}$ is a subset of a smaller subset of $\{1,\ldots,n\}$.

\medskip

Our first basic result specifies much further the nature of the minimal binomial generators.

\begin{Theorem}\label{format}
Let $I\subset R=k[x_1,\ldots,x_n]$ be a uniform monomial ideal as above.
Then the polynomial ring $S:=R[y_1,\ldots,y_n,w]$ admits a grading under which the presentation ideal $\mathcal{I}$ of the Rees algebra of $I$ over it is generated by homogeneous binomials in this grading.

Moreover:
\begin{itemize}
\item[{\rm (a)}] If $a\leq nb$, letting $1\leq p\leq n$ be the unique integer such that $(p-1)b<a\leq pb$, then
any minimal binomial generator of external degree $\delta$ can be written in the form
\begin{equation}
(x_{i_1}\cdots x_{i_{\delta}})^{a-\delta b}\,w^{\delta}- P(i_1,\ldots,i_{\delta})^{\delta b}\,y_{i_1}\cdots y_{i_{\delta}},
\end{equation}
where  $\delta\leq p$,
 with the convention that if $\delta=p$ then the $\xx$-term on the left hand side goes over to the right hand side with exponent $-(a-\delta b)=\delta b-a$.
\item[{\rm (b)}] If $a>nb$, then
 any minimal binomial generator of external degree $\delta$ can be written in the form
\begin{equation}
(x_{i_1}\cdots x_{i_{\delta}})^{a-\delta b}\,w^{\delta}- P(i_1,\ldots,i_{\delta})^{\delta b}\,y_{i_1}\cdots y_{i_{\delta}},
\end{equation}
where $\delta\leq n$.
{\rm (}no convention needed in this case since for $\delta=n$, there is no $\xx$-term on the right hand side{\rm )}.
\end{itemize}
\end{Theorem}
\demo
Start with generators of the presentation ideal of the symmetric algebra of $I$.
It is easy to see that the syzygies of $I$ are generated by the Koszul relations of the pure powers $x_1^a,\ldots, x_n^a$ and by the reduced relations of $(x_1\cdots x_n)^b$ with each one of the pure powers.
In other words,  ${\mathcal L}\subset S=R[y_{1},\ldots,y_{n},w]$ is generated by the binomials $$K_{i,j}=x_{i}^{a}y_{j}-x_{j}^{a}y_{i},~i,j\in\{1,\ldots,n\},~i<j,$$ $$L_{i}=x_i^{a-b}w-P(i)^{b}y_{i},~i\in \{1,\ldots,n\}.$$
Now, these binomials are homogeneous in $S$ by
attributing the following weights to the variables: $\deg(x_i)=1$ and
$\deg(w)=nb-a+1,\deg(y_j)=1$ if $a\leq nb$, while $\deg(w)=1, \deg(y_j)=a-nb+1$ if $a\geq nb$.
Therefore, ${\mathcal L}$ is homogeneous for these weights.
Since ${\mathcal I}={\mathcal L}:I^{\infty}$ and $I$ is monomial, it follows that ${\mathcal I}$ is generated by binomials which are homogeneous as well under the same weights.
Indeed, one has the string of inclusions
$${\mathcal I}={\mathcal L}:I^{\infty}\subset {\mathcal L}:(x_1)^{\infty}\subset {\mathcal I}:(x_1)^{\infty}={\mathcal I},$$
the last equality because ${\mathcal I}$ is a prime ideal.
Then by \cite[Corollary 1.7 (a)]{binomials} (or, directly, by \cite[Corollary 1.9]{binomials}), $\mathcal I$ is generated by binomials and hence by homogeneous binomials as $x_1$ is homogeneous of degree $1$.
(Note that the counterexamples in \cite{binomials} are non-prime.)

By Lemma~\ref{basic_form},
a binomial in $\mathcal I$ belonging to a set of minimal generators thereof is of the form
$$\mathbf{m}(\xx)w^{\delta}-\mathbf{n}(\xx)y_{i_1}^{\alpha_{i_1}}\cdots y_{i_s}^{\alpha_{i_s}},$$
with $1\leq i_1<\cdots < i_s\leq n$, $\alpha_{i_j}>0$, and $\mathbf{m}(\xx), \mathbf{n}(\xx)$ suitable monomials in $R$ such that $\gcd \{\mathbf{m}(\xx), \mathbf{n}(\xx)\}=1$.

\smallskip

In addition, one has the following three basic principles:

\smallskip

$\bullet$  $w$ corresponds to a monomial that involves all variables of $R$; this implies that the monomial $\mathbf{n}(\xx)$ must involve the variables indexed by the complementary subset $\{j_{s+1},\ldots, j_n\}:=\{1,\ldots,n\}\setminus \{i_1,\ldots , i_s\}$ and, since $\gcd \{\mathbf{m}(\xx), \mathbf{n}(\xx)\}=1$, the variables effectively involved in $\mathbf{m}(\xx)$ must be indexed by a subset of $\{i_1,\ldots , i_s\}$.
Therefore, the monomial has the form
$$x_{i_1}^{d_{i_1}}\cdots x_{i_s}^{d_{i_s}}\,w^{\delta} - x_{i_{s+1}}^{c_{i_{s+1}}}\cdots x_{i_n}^{c_{i_n}}y_{i_1}^{\alpha_{i_1}}\cdots y_{i_s}^{\alpha_{i_s}}$$
for suitable  exponents $d_{i_l}\geq 0$, for $l=1,\ldots,s$ (some of which may vanish) and $c_{i_k}$, for $k=s+1,\ldots,n$ (which are positive).

\smallskip

$\bullet$ Weighted homogeneity implies the equalities
\begin{equation}\label{quasi-homog_case1}
(nb-a+1)\delta + \sum_{l=1}^s d_{i_l} = \sum_{l=1}^s \alpha_{i_l} +\sum_{k=s+1}^n c_{i_k}
\end{equation}
if $a\leq nb$, and
\begin{equation}\label{quasi-homog_case2}
\delta + \sum_{l=1}^s d_{i_l} = (a-nb+1)\sum_{l=1}^s \alpha_{i_l} +\sum_{k=s+1}^n c_{i_k}
\end{equation}
if $a\geq  nb$.

Moreover, since upon evaluation the powers $x_{i_k}^{c_{i_k}}$ on the right hand side can only cancel against the ones coming from $w^{\delta}$ on the left hand side, we see that ${c_{i_k}}=\delta b$ for every $k=s+1,\ldots,n$.
By the same token, $d_{i_l}=a \alpha_{i_l}-\delta b$ for every $l=1,\ldots,s$.

\smallskip

$\bullet$ Lastly, since the Rees algebra ${\mathcal R}_R(I)$ is also standard graded over $R={\mathcal R}_R(I)_0$, we may assume that the binomial is homogeneous with respect to the external variables (however, we warn that ${\mathcal R}_R(I)$ is {\em standard}  bigraded over $k$ if and only if $a=nb$).
This means that $\delta= \sum_{l=1}^s \alpha_{i_l}$, a formula already found in the above lemma.

\smallskip

So we can assume our binomial to look like $$x_1^{a\alpha_1-\delta b}\cdots x_s^{a\alpha_s-\delta b}w^{\delta}-(x_{s+1}\cdots x_n)^{\delta b}y_1^{\alpha_1}\cdots y_s^{\alpha_s}, \alpha_i\geq 1.$$

\medskip

Suppose $\delta\geq p+1$. The goal is to show that this binomial can be generated by binomials in $\mathcal I$ with $w$ raised to a power $\leq p$. Since $a<(p+1)b$ and $a\alpha_i-\delta b>0$, then $\alpha_i\geq 2$ for all $i=1,\ldots,s$.

\smallskip

If $s\geq p$,  consider the polynomial
$$H:=w^p-(x_1\cdots x_p)^{pb-a}(x_{p+1}\cdots x_n)^{pb}y_1\cdots y_p\in\mathcal I.$$
If $a=pb$, consider $H:=w^p-y_1\cdots y_p$. By primality of $\mathcal I$, using $H$, our binomial is generated by $H$ and by the following binomial in $\mathcal I$
$$x_1^{a(\alpha_1-1)-(\delta-p)b}\cdots x_p^{a(\alpha_p-1)-(\delta-p)b}x_{p+1}^{a\alpha_{p+1}-(\delta-p)b} \cdots x_s^{a\alpha_s-(\delta-p)b}w^{\delta-p}$$ $$-(x_{s+1}\cdots x_n)^{(\delta-p)b}y_1^{\alpha_1-1}\cdots y_p^{\alpha_p-1}y_{p+1}^{\alpha_{p+1}}\cdots y_s^{\alpha_s},$$
where $w$ is raised to $\delta-p$, and in addition the exponents of $x_i$ on the left do not vanish  since $a\alpha_i>\delta b$, then $a(\alpha_i-1)-(\delta-p)b>pb-a\geq 0$.

\smallskip

If $s\leq p-1$,  consider
$$G:=(x_1\cdots x_s)^{a-sb}w^s-(x_{s+1}\cdots x_n)^{sb}y_1\cdots y_s.$$
Then, by the same token as above, using $G$, the binomial can be generated by $G$ and by the following  binomial in $\mathcal I$: $$x_1^{a(\alpha_1-1)-(\delta-s)b}\cdots x_s^{a(\alpha_s-1)-(\delta-s)b}w^{\delta-s}- (x_{s+1}\cdots x_n)^{(\delta-s)b}y_1^{\alpha_1-1}\cdots y_s^{\alpha_s-1}.$$

\smallskip

Recursively, in both situations above ($s\geq p$ and $s\leq p-1$), our binomial can be generated by binomials in $\mathcal I$ of the same shape with  $w$ raised to a power $\leq p$.

\smallskip

The concluding blow is given by the following result:

{\sc Claim.} With the preceding notation, if $\delta\leq p$, then we can assume $\alpha_1=\cdots=\alpha_s=1$, and $s=\delta$.

For the proof, assume $\alpha_1\geq 2$. Then $a\alpha_1-\delta b\geq 2a-\delta b=a-b +a-(\delta-1)b$. Since $p\geq \delta$ and $a>(p-1)b$, then $a-(\delta-1)b>0$. Our binomial can be written as

$$x_1^{a(\alpha_1-1)-(\delta-1)b}x_2^{a\alpha_2-\delta b}\cdots x_s^{a\alpha_s-\delta b}w^{\delta-1}\underbrace{(L_1+(x_2\cdots x_n)^by_1)}_{x_1^{a-b}w}-(x_{s+1}\cdots x_n)^{\delta b}y_1^{\alpha_1}\cdots y_s^{\alpha_s}.$$

Since $L_1\in\mathcal I$ and we only care for minimal generators, by simplifying by $(x_{s+1}\cdots x_n)^by_1$ one can assume the binomial to be of the form $$x_1^{a(\alpha_1-1)-(\delta-1)b}x_2^{a\alpha_2-(\delta-1)b}\cdots x_s^{a\alpha_s-(\delta-1) b}w^{\delta-1}-(x_{s+1}\cdots x_n)^{(\delta-1)b}y_1^{\alpha_1-1}y_2^{\alpha_2}\cdots y_s^{\alpha_s},$$
where both $\alpha_1$ and  $\delta$ dropped by $1$.
Therefore, recursion takes care of the conclusion.

The case where $a>nb$ is handled similarly.

This concludes the proof of the claim and also of the theorem.
\qed

\subsection{Sylvester forms as generators}\label{generation}

For the reader's convenience, we recall once more the following notation: if $\{i_1,\ldots, i_j\}$ is a subset of $\{1,\ldots,n\}$ in the natural order of the integers, we denote by $P(i_1,\ldots,i_j)$ the product of the variables in the complementary set $\{x_1,\ldots,x_n\}\setminus \{x_{i_1},\ldots, x_{i_j}\}$.

The next theorem partly summarizes the results of the preceding part, adding information on the nature of the generators as Sylvester forms.

\begin{Theorem}\label{gens_relation_type}
Let $I\subset R$ be a uniform monomial ideal as above and let $r$ denote its reduction number as established in {\rm Proposition~\ref{red_numbers}}.
Then:

\smallskip

{\rm (a)}  $\mathcal{I}$ is  generated by
$$ {{n}\choose {2}}+ \sum_{\delta=1}^r  {{n}\choose {\delta}} +1,
$$
quasi-homogeneous binomials, where $r$ is the reduction number of $I${\rm ;} of these, $ {{n}\choose {2}}$ are the Koszul syzygies of the generators of $I$ and the remaining ones are each a binomial of the form
$$(x_{i_1}\cdots x_{i_{\delta}})^{a-\delta b}\,w^{\delta}- P(i_1,\ldots,i_{\delta})^{\delta b}\,y_{i_1}\cdots y_{i_{\delta}},
$$
where $1\leq \delta\leq r+1$ {\rm (}with the same convention as stated in {\rm Theorem~\ref{format}} in the case $a\leq nb${\rm )}.

{\rm (b)} Moreover, each binomial in the previous item is a Sylvester form obtained in an iterative form out of the syzygy forms.

\smallskip

{\rm (c)} The relation type of $I$ is $r+1$.

\end{Theorem}
\demo
(a) The proof of the generation statement will consist in showing that a quasi-homogeneous generator of $\mathcal I$ of arbitrary standard degree in the external variables $y_1,\ldots,y_n, w$ belongs to the ideal generated by the binomials in the statement, with standard external degrees bounded by the reduction number of $I$.
Thus, the result will be a consequence of Theorem~\ref{format} and of Proposition~\ref{red_numbers}.

\smallskip

From the above degree reduction result and from Theorem~\ref{format} we deduce that, for each $2\leq \delta \leq r$, where $r$ is the reduction number of $I$, $\mathcal{I}$ admits
${{n}\choose {\delta}}$
generators which are quasi-homogeneous binomials.
Generators for $\delta=1$ are the syzygy binomials, which add up ${{n}\choose {2}}+ n$ generators in standard degree $1$.

Finally, we deal with generators in standard degree $r+1$.
In the case where $a> nb$, then there is a unique generator in degree $n$ given in Theorem~\ref{format}, namely, $(x_1\cdots x_n)^{a-nb}w^n-y_1\cdots y_n$.
In the case where $a\leq nb$ and $p\leq n$ is the unique integer such that $(p-1)b<a\leq pb$, we obtain ${{n}\choose {p}}$ generators, one for each choice of an ordered subset $\{i_1,\ldots, i_p\}\subset \{1,\ldots,n\}$:
$$S_{i_1,\ldots, i_p}:=w^p-P(i_1,\ldots,i_p)^{pb}\, (x_{i_1}\cdots x_{i_p})^{pb-a}\, y_{i_1}\cdots y_{i_p}.
$$
We now show that fixing one of these, the remaining ones belong to the ideal generated by this one and the Koszul relations.
To prove this assertion it suffices to fix one subset $\{i_1,\ldots, i_p\}$ and another subset obtained by one transposition.
Without loss of generality, we assume the fixed subset is $\{1,\ldots, p\}$ and the other one is $\{1,\ldots, p-1,p+1\}$.

\smallskip

{\sc Claim:} With the above notation and the previous notation for the Koszul relations, one has
$$S_{i_1,\ldots, i_{p-1}, i_{p+1}}= S_{i_1,\ldots, i_p} + M(\xx)\,y_2\cdots y_p\, L_{1,p+1} - M(\xx)\,y_2\cdots y_{p-1} y_{p+1}\, L_{1,p},
$$
where $M(\xx)=(x_1\cdots x_px_{p+1})^{pb-a}(x_{p+2}\cdots x_n)^{pb}$.

The proof is a straightforward calculation by developing the right hand side.

As a consequence, also for the case $(p-1)b<a\leq pb$ there is a unique minimal generator in standard degree $p$.
Summing up, in both cases, we get
$$ {{n}\choose {2}}+ \sum_{\delta=1}^r  {{n}\choose {\delta}} +1
$$
minimal quasi-homogeneous binomial generators.

\medskip

(b) We next show that the generators of the first part are indeed Sylvester forms obtained iteratively.

Recall once more the form of the generators of
${\mathcal L}\subset S=R[y_{1},\ldots,y_{n},w]$: the Koszul relations \begin{equation}\label{koszul}
K_{i,j}=x_{i}^{a}y_{j}-x_{j}^{a}y_{i},~i,j\in\{1,\ldots,n\},~i<j
\end{equation}
and the reduced (Taylor) relations
\begin{equation}\label{taylor}
L_{i}=x_i^{a-b}w-P(i)^{b}y_{i},~i\in \{1,\ldots,n\}.
\end{equation}

We start by availing ourselves of Sylvester forms of degree $2$.
For this, take any two distinct indices $l,i\in \{1,\ldots,n\}$, say, $l<i$. We form the Sylvester content matrix of $\{L_{l},L_{i}\}$ with respect to the complete intersection $\{x_{l}^{b},x_{i}^{b}\}$:
$$\left[
\begin{array}{c}
L_{l}  \\ [5pt]
L_{i}
\end{array}
\right] = \left[
\begin{array}{c}
 x_{l}^{a-b}w-P(l)^{b}y_{l}  \\ [5pt]
 x_{i}^{a-b}w-P(i)^{b}y_i
\end{array}
\right]=\underbrace{\left(
\begin{array}{cc}
x_{l}^{a-2b}w & -P(l,i)^{b}y_{l}   \\ [5pt]
 -P(l,i)^{b}y_{i} & x_{i}^{a-2b}w
\end{array}
\right)}_{M_{2}^{l,i}}\left[
\begin{array}{c}
x_{l}^{b}  \\ [5pt]
x_{i}^{b}
\end{array}
\right].$$
Set $H_{2}^{l,i}=\mathrm{det}(M_{2}^{l,i})=(x_{l}x_{i})^{a-2b}w^{2}-P(l,i)^{2b}y_{l}y_{i}.$
Note that, since we are assuming that $a>2b$, we obtain this way ${{n}\choose {2}}$ distinct forms of external degree $2$.

We now induct on the degree. Thus, suppose that for $j\in\{1,\ldots,n\}$ with $a> jb$, one has found ${{n}\choose {j}}$
Sylvester forms, of external degree $j$, each of the shape \begin{equation}\label{general_sylvester}
H_{j}^{i_{1},\ldots,i_{j}}=(x_{i_{1}}x_{i_{2}}\cdots x_{i_{j}})^{a-jb}w^{j}-P(i_1,\ldots,i_j)^{jb}y_{i_{1}}\cdots y_{i_{j}},
\end{equation}
with $i_{1},\ldots,i_{j}\in \{1,\ldots,n\}$ and $i_{1}<\cdots<i_{j}.$
Then for every $l\in \{1,\ldots,n\}\setminus \{i_{1},\ldots,i_{j}\}$, we obtain a Sylvester content matrix of $L_l, H_{j}^{i_{1},\ldots,i_{j}}$ with respect to the complete intersection $(x_{l}^{jb},(x_{i_{1}}\cdots x_{i_{j}})^{b})$:
{\small
\begin{eqnarray}\nonumber
\left[
\begin{array}{c}
L_{l}  \\ [5pt]
H_{j}^{i_{1},\ldots,i_{j}}
\end{array}
\right] &=& \left[
\begin{array}{c}
 x_{l}^{a-b}w-P(l)^{b}y_{l}  \\ [5pt]
 (x_{i_{1}}\cdots x_{i_{j}})^{a-jb}w^{j}-P(i_1,\ldots,i_j)^{jb}y_{i_{1}}\cdots y_{i_{j}}
\end{array}
\right]\\ [8pt] \nonumber
&=&\underbrace{\left(
\begin{array}{cc}
x_{l}^{a-(j+1)b}w & -P(i_1,\ldots,i_j,l)^{b}y_{l}   \\ [5pt]
 -P(i_1,\ldots,i_j,l)^{jb}y_{i_{1}}\cdots y_{i_{j}} & (x_{i_{1}}\cdots x_{i_{j}})^{a-(j+1)b}w^{j}
\end{array}
\right)}_{M_{j+1}^{i_{1},\ldots,l,\ldots,i_{j}}}\left[
\begin{array}{c}
x_{l}^{jb}  \\ [5pt]
(x_{i_{1}}\cdots x_{i_{j}})^{b}
\end{array}
\right].
\end{eqnarray}
}

This yields a new Sylvester form of external degree $j+1$: $H_{j+1}^{i_{1},\ldots,\ldots,i_{j},\,l}=\mathrm{det}(M_{j+1}^{i_{1}, \ldots,\ldots,i_{j},\,l})$ $$\hspace{1,8cm}=(x_{i_{1}} \cdots x_{i_{j}},\,x_{l})^{a-(j+1)b}w^{j+1}-P(i_1,\ldots,i_j,l))^{(j+1)b}y_{i_{1}}\cdots y_{i_{j}}\cdots y_l.$$
(Here, we assume that $\{i_1,\ldots,i_j,l\}$ is written in increasing order.)
This way we have produced ${{n}\choose {j+1}}$ distinct Sylvester forms of external degree $j+1$.

To conclude the inductive procedure, we divide the proof into the two basic cases:

\smallskip

{\bf (i)} $a\leq nb$.

\smallskip

In this case, let $1\leq p\leq n$ be the smallest integer such that $(p-1)b<a\leq pb.$
By the previous argument, since $a>(p-1)b$ then a  Sylvester form of standard degree $(p-1)$ over $R$ has the shape $$H_{p-1}^{i_{1},\ldots,i_{p-1}}=(x_{i_{1}}\cdots x_{i_{p-1}})^{a-(p-1)b}w^{p-1}-P(i_1,\ldots,i_{p-1})^{(p-1)b}y_{i_{1}}\cdots y_{i_{p-1}}, $$
with $\{i_{1},\ldots,i_{p-1}\}$ an ordered subset of $\{1,\ldots,n\}$.   Take the Sylvester form of $\{L_{l},H_{p-1}^{i_{1},\ldots,i_{p-1}}\}$ with respect to $\{ x_{l}^{a-b},
(x_{i_{1}}\cdots x_{i_{p-1}})^{a-(p-1)b} \}$, since $a\leq pb$:

$$
\left[
\begin{array}{c}
L_{l}  \\
H_{p-1}^{i_{1},\ldots,i_{p-1}} \\
\end{array}
\right]= M_{p}^{i_{1},\ldots,l,\ldots,i_{p-1}}
\cdot \left[
\begin{array}{c}
x_{l}^{a-b}  \\
(x_{i_{1}}\cdots x_{i_{p-1}})^{a-(p-1)b}  \\
\end{array}
\right],
$$
where $ M_{p}^{i_{1},\ldots,l,\ldots,i_{p-1}}$ denotes the content matrix
{\small
$$
\left(
\begin{array}{cc}
w & -P(i_1,\ldots, i_{p-1},l)^{b}(x_{i_{1}}\cdots x_{i_{p-1}})^{pb-a}y_{l}   \\ [8pt]
 -P(i_1,\ldots, i_{p-1},l)^{(p-1)b}x_{l}^{pb-a}y_{i_{1}}\cdots y_{i_{p-1}} & w^{p-1} \\
\end{array}
\right).
$$
}

Thus,
\begin{eqnarray}\nonumber
H_{p}^{i_{1},\ldots,l,\ldots, i_{p-1}}&=&\mathrm{det}(M_{p}^{i_{1},\ldots,l,\ldots,i_{p-1}})\\ \nonumber
&=& w^{p}-P(i_1,\ldots, i_{p-1},l)^{pb}(x_{i_{1}}\ldots x_{l}\ldots x_{i_{p-1}} )^{pb-a}y_{i_{1}}\cdots y_{l}\cdots y_{i_{p-1}}.
\end{eqnarray}

{\bf (ii)} $a> nb.$

\smallskip

By the previous argument, since $a> nb$ then a Sylvester form of standard degree $n$ over $R$ has the shape $$H_{n}^{1,\ldots,n}=(x_{1}\cdots x_{n})^{a-nb}w^{n}-y_{1}\cdots y_{n}.$$

\medskip

(c) This follows immediately from the details of the generation as described in (a).
\qed

\begin{Remark}\rm
Note the sharp difference between cases (i) and (ii) at the end of the proof above: if $p=n$ then there is a unique binomial Sylvester form with a term a pure power of $w$ (namely, $w^n$), while for $p<n$ there are various such binomials having the pure term $w^p$ -- although only one emerges as part of a minimal set of generators, as explained in the proof of the previous theorem.
\end{Remark}

\section{Combinatorial structure of the Rees ideal}

We keep the notation of the previous part.
Recall that, given an integer $2\leq j\leq p-1$, where $p-1\leq n-1$ is the reduction number of the ideal $I\subset S=k[x_1,\ldots,x_n]$, and an increasing sequence of integers $i_{1}<\cdots<i_{j}$ in $\{1,\ldots,n\}$, we had a well-defined Sylvester form $H_{j}^{i_{1},\ldots,i_{j}}$ in the set of generators of the Rees ideal $\mathcal{R}_R(I)$.
This polynomial is weighted homogeneous in all concerned variables and homogeneous of degree $j$ in the presentation variables $y_1,\ldots,y_n, w$.
We will order the set of these forms in the following way: first, if two of these forms $H_{j}^{i_{1},\ldots,i_{j}}$ and  $H_{j}^{k_{1},\ldots,k_{j}}$ have the same presentation degree $j$ then we set $H_{j}^{i_{1},\ldots,i_{j}}$ before $H_{j}^{k_{1},\ldots,k_{j}}$ provided $i_r<k_r$, where $r$ is the first index from the left
such that $i_r\neq k_r$; second, we decree that the last form $H_{j}^{n-j+1,\ldots,n}$ of degree $j$ in this ordering precedes the first form $H_{j}^{1,2,\ldots,j+1}$ of the next presentation degree $j+1$.

The presentation ideal of the symmetric algebra of $I$ is denoted $\mathcal{L}$ as before. It is generated by the Koszul relations $K_{i,j},\, 1\leq i<j\leq n$ and the reduced Taylor relations $L_i,\, 1\leq i\leq n$, as in (\ref{koszul}) and (\ref{taylor}).

\smallskip
We will need the following easy properties of the colon ideal in the proof of the next proposition:
\begin{Lemma}\label{lemma_colon}
Let $J\subset R$ be an ideal in a ring and $f\in R$.
Then:
\begin{enumerate}
\item[{\rm (i)}] $(J:f)f=J\cap (f)$.
\item[{\rm (ii)}] Suppose that $R$ is a polynomial ring over a field and $<$ is a monomial order. Then $\mathrm{in}_<(J:f)\subset \mathrm{in}_<(J):\mathrm{in}_<(f)${\rm ;} if in addition $\mathrm{in}_<(J):\mathrm{in}_<(f)\subset J:f$ then the equality $\mathrm{in}_<(J):\mathrm{in}_<(f)= J:f$ holds.
\end{enumerate}
\end{Lemma}
\demo (i) This is straightforward from the definition of the colon ideal.

(ii) The inclusion $\mathrm{in}_<(J:f)\subset \mathrm{in}_<(J):\mathrm{in}_<(f)$ follows immediately from the definition of the initial ideal.

Now let $F\in J:f$. Then, by the above inclusion and the assumption, one has $\mathrm{in}_<(F)\in J:f$, hence $G:=F-\mathrm{in}_<(F)\in J:f$.
By induction on the number of nonzero terms of a polynomial in $R$, we have $G\in \mathrm{in}_<(J):\mathrm{in}_<(f)$. It follows that $F\in \mathrm{in}_<(J):\mathrm{in}_<(f)$.
\qed

\subsection{Initial ideals}

In the following propositions we discuss the preliminaries on Gr\"obner basis and initial ideals related to the ordered sequence of Sylvester forms.

\begin{Proposition}\label{G-basis}
Let  $2\leq j\leq p-1$, where $p-1\leq n-1$ is the reduction number of the ideal $I\subset S=k[x_1,\ldots,x_n]$, and let $i_{1}<\cdots<i_{j}$ be an ordered subset of  $\{1,\ldots,n\}$.
The set
$$\Sigma(i_{1},\ldots,i_{j}):= \{K_{i,k} \,(1\leq i<k\leq n), L_i\, (1\leq i \leq n),H_{2}^{1,2},\ldots,H_{j}^{i_{1},\ldots,i_{j}}\}
$$
is a Gr\"obner basis of the ideal ${\mathcal H}(i_{1},\ldots,i_{j}):=(\mathcal{L},H_{2}^{1,2},\ldots,H_{j}^{i_{1},\ldots,i_{j}})$
in the lexicographic order on $w>x_n>\cdots >x_1>>\cdots$. In particular, the initial ideal of ${\mathcal H}(i_{1},\ldots,i_{j})$ is
\begin{equation}\label{initial_ideal}
\{x_{i}^{a}y_{k}\, (1\leq i<k\leq n),\,x_{i}^{a-b}w\, (1\leq i\leq n),\, (x_{1}x_{2})^{a-2b}w^{2},\ldots,(x_{i_{1}}\cdots x_{i_{j}})^{a-jb}w^{j}\},
\end{equation}
where $j$ and $\{x_{i_1},\ldots, x_{i_j}\}$ flow as in the statement.
\end{Proposition}
Since the proof is a case-by-case scrutiny of $S$-pairs, we postpone it to the last section of the paper.

\begin{Proposition}\label{colon-of-initials} With the above setting, let $H_{j'}^{k_1,\cdots , k_{j'}}\in S$ denote the first Sylvester form succeeding the Sylvester from $H_{j}^{i_1,\cdots, i_j}\in S$ in the prescribed ordering of these forms.
\begin{enumerate}
\item[{\rm (a)}]
If $j=j'$, one has
{\small
\begin{eqnarray} \label{initial_colon1}\nonumber
\lefteqn { {\rm in}({\mathcal H}(i_{1},\ldots,i_{j})):  {\rm in} (H_{j}^{k_{1},\ldots,k_{j}})=\left(x_{k_{1}}^{(j-1)b},\ldots,x_{k_{j}}^{(j-1)b}, x_{u}^{a-jb}, x_{k_{j}+1}^{a-b},\ldots,x_{n}^{a-b},
(x_{r_{1}}\cdots x_{r_{s}})^{(j-s)b},\right.}\kern7cm\\ \nonumber
&& \left.(x_{q_{1}}\cdots x_{q_{r}})^{(j-s)b}(x_{d_{1}}\cdots x_{d_{s-r}})^{a-sb} \right)S,
\end{eqnarray}
}
for all $k_i<u<k_{i+1},i=1,\ldots,j-1$ and all choices of indices $s\in\{1,\ldots,j-1\},~r\in\{0,\ldots,s-1\}$,
 of ordered subsets
$\{r_{1}<\cdots<r_{s}\}\subset \{k_1,\ldots,k_j\}$, $\{q_{1}<\cdots<q_{r}\}\subset \{k_1,\ldots,k_j\}$  and of an ordered set $d_{1}<\cdots<d_{s-r}$ with $k_j<d_{1}$.
\item[{\rm (b)}]
If $j'=j+1$ {\rm (}and hence $\{k_1,\cdots , k_{j'}\}=\{1,\ldots,j+1\}${\rm )}, one has
{\small
\begin{eqnarray} \label{initial_colon2}\nonumber
\lefteqn {\kern-2.5cm{\rm in}({\mathcal H}(i_{1},\ldots,i_{j})):{\rm in}(H_{j+1}^{1, \ldots,j+1}) =\left( x_{1}^{jb}, \ldots,x_{j+1}^{jb},x_{j+2}^{a-b},\ldots,x_{n}^{a-b},(x_{r_{1}}\cdots x_{r_{s}})^{(j+1-s)b},\right.}\kern0.5cm\\ \nonumber
&&\left. (x_{q_{1}}\cdots x_{q_{r}})^{(j+1-s)b}(x_{d_{1}}\cdots x_{d_{s-r}})^{a-sb}\right)S,
\end{eqnarray}
}
for all choices of indices $s\in\{1,\ldots,j\},~r\in\{0,\ldots,s-1\}$,
 of ordered subsets
$\{r_{1}<\cdots<r_{s}\}\subset \{1,\ldots,j+1\}$, $\{q_{1}<\cdots<q_{r}\}\subset \{1,\ldots,j+1\}$  and of an ordered set $d_{1}<\cdots<d_{s-r}$ with $j+1<d_{1}$.

\smallskip

{\rm (}In both cases, we adopt the convention that  $x_{q_{0}}=1$.{\rm )}
\end{enumerate}
\end{Proposition}

For both items, we will apply Lemma~\ref{lemma_colon} (i), by which one is to compute a minimal set of generators of the intersection of the two initial ideals on the left hand side, then divide each generator by the initial term of $H_{j'}^{k_{1},\ldots,k_{j'}}$.
To get a minimal set of generators of the intersection we use a well-known principle  (see, e.g., \cite[Proposition 1.2.1]{HH}), by which this set is the set of the least common multiples of ${\rm in}(H_{j'}^{k_{1},\ldots,k_{j'}})$ and each minimal generator of ${\rm in}({\mathcal H}(i_{1},\ldots,i_{j}))$.

The details of the proof are given in the last section.


\medskip

The next result is slightly surprising.

\begin{Proposition}\label{successive_colon}
With the previously established notation, one has
$$\mathcal{H}(i_{1}, \ldots,i_{j}):H_{j'}^{k_{1},\ldots,k_{j'}}={\rm in}(\mathcal{H}(i_{1}, \ldots,i_{j})):{\rm in}(H_{j'}^{k_{1}, \ldots,k_{j'}}).$$
In particular, the colon ideal on the left hand side is a monomial ideal.
\end{Proposition}

The proof hinges on the explicit form of the generators given in the previous proposition.
The computation is again a case-by-case calculation and quite often it requires some ingenuity as to how the generator looks and how the result of the calculation ought to look like.
Since at this point it will give no additional conceptual contribution to the rhythm of the exposition, we postpone the details to last section.

\subsection{Almost Cohen--Macaulayness}
In this part we deal with the depth of the Rees algebra of the ideal $I\subset S$, which wraps-up the main goal of the paper.

 In the notation  of the preceding sections, the main result is as follows.

\begin{Theorem}\label{rees_is_acm}
Let $I\subset R=k[x_1,\ldots,x_n]$ denote a uniform monomial ideal as in {\rm Section~\ref{generation}}. Then $S/\mathcal{H}(i_1,\ldots,i_j)$ has depth at least $n$ for every tuple $i_1<\cdots<i_j$. In particular, the Rees algebra $\mathcal{R}_R(I)$ of $I$ is an almost Cohen--Macaulay ring.
\end{Theorem}
\demo We basically follow the idea of \cite[Theorem  3.14 (b)]{sylStefan}.
Namely, produce a sequence of mapping cones, each a free resolution of the sequential ideal $$\mathcal{H}(i_1,\ldots,i_j):=(\mathcal{L},H_{2}^{1,2},\ldots,H_{j}^{i_1, \ldots,i_j})$$
discussed above, ending with a free resolution of $\mathcal{R}_R(I)$; at each step the mapping cone has length at most $n+1$. Therefore, the depth of $\mathcal{R}_R(I)$ will turn out to be at least $2n+1-(n+1)=n$, as desired.

In a precise way, we now argue that for each tuple $i_1<\cdots<i_j$, starting from the first tuple $1<2$, a free $S$-resolution of $S/\mathcal{H}(k_1, \ldots,k_{j'})$ is the mapping cone of the map of complexes from a resolution of $S/(\mathcal{H}(i_1,\ldots,i_j):H_{j'}^{k_1, \ldots,k_{j'}})$ to a resolution of $S/\mathcal{H}(i_1,\ldots,i_j)$ induced by multiplication by $H_{j'}^{k_1, \ldots,k_{j'}}$ on $S$, where $k_1< \cdots <k_{j'}$ is the first tuple succeeding $i_1<\cdots <i_j$ in the ordering explained before.

To see this, we induct on the number of generators $\mathcal{H}(i_1,\ldots,i_j)$.

Now, by Proposition~\ref{successive_colon} and Proposition~\ref{colon-of-initials}, the generators of the colon ideal $\mathcal{H}(i_1,\ldots,i_j):H_{j'}^{k_1, \ldots,k_{j'}}$ are elements of $R$ containing powers of all variables. Therefore, these monomials generate an $R_+$-primary ideal of $R$, and hence a free $S$-resolution of $S/(\mathcal{H}(i_1,\ldots,i_j):H_{j'}^{k_1, \ldots,k_{j'}})$ is obtained by flat base change $R\subset S$ from a minimal free $R$-resolution of length $n$.

In the first step one has $\mathcal{H}(1,2)=(\mathcal{L},H_{2}^{1,2})$.
Since the ideal $I\subset R$ is an almost complete intersection of finite length, $S/\mathcal{L}$ is Cohen--Macaulay (\cite[Corollary 10.2]{Trento}). As the codimension of the Rees algebra of $I$ on $S$ is $n$, the codimension of  $S/\mathcal{L}$ is at least $n$. But since $\mathcal{L}\subset R_+=(x_1,\ldots,x_n)S$ then the codimension is $n$.

We consider the map of complexes induced by multiplication by $H_{2}^{1,2}$ on $S$:
$$\begin{array}{ccccccccccccc}
0 & \rar & S^{\beta_n}  & \lar &\cdots &\lar & S^{\beta_2} & \lar & S^{\beta_1} & \lar & S & \rar & 0\\
&&\uparrow &&&& \uparrow && \uparrow && \uparrow && \\
0 & \rar & S^{\alpha_n}  & \lar &\cdots &\lar & S^{\alpha_2} & \lar & S^{\alpha_1} & \lar & S & \rar & 0
\end{array},
$$
where the upper complex is a free resolution of $S/\mathcal{L}$ and the lower one  is the free $S$-resolution of $S/(\mathcal{L}:H_{2}^{1,2})$ extended from the free $R$-resolution by flat base change $R\subset S$.
(Note that $\beta_1={{n+1}\choose {2}}$ is the minimal number of generators of $\mathcal{L}$, but all the remaining Betti number of both resolutions are harder to guess.)

The mapping cone is a free $S$-resolution of $S/\mathcal{H}(1,2)$ (not minimal as there will be cancellation in general).
By definition, this $S$-resolution has length at most $n+1$.

The general step of the induction is entirely similar, by taking the mapping cone of the map of complexes induced by multiplication by $H_{j'}^{k_1, \ldots,k_{j'}}$:
$$\begin{array}{ccccccccccccccc}
0 & \rar & S^{\beta_{n+1}} & \lar & S^{\beta_n}  & \lar &\cdots &\lar & S^{\beta_2} & \lar & S^{\beta_1} & \lar & S & \rar & 0\\
&&\uparrow && \uparrow &&&& \uparrow && \uparrow && \\
&& 0 & \lar & S^{\alpha_n}  & \lar &\cdots &\lar & S^{\alpha_2} & \lar & S^{\alpha_1} & \lar & S & \rar & 0
\end{array},
$$
where the upper complex is a (not necessarily minimal) free resolution of $S/\mathcal{H}(i_1,\ldots,i_j)$ and the lower one is the $S$-resolution of $S/(\mathcal{H}(i_1,\ldots,i_j):H_{j'}^{k_1, \ldots,k_{j'}})$ extended by flat base change from  a minimal free $R$-resolution.
Here we have used for simplicity the same notation for the Betti number as above, but of course they are different.

Because the lower complex has length at most the length of the upper complex, the mapping cone is again a free $S$-resolution of length at most $n+1$.

By Theorem~\ref{gens_relation_type} and the previous discussion of this section, the presentation ideal $\mathcal{I}$ of the Rees algebra on $S$ is the sequential ideal $\mathcal{H}(1,\ldots,p)$, where $p-1$ is the reduction number of $I$. Therefore the above gives that $\mathcal{R}_R(I)$ has an $S$-resolution of length at most $n+1$, as was to be shown.
\qed

\section{Proofs}

\subsection{Proof of Proposition~\ref{G-basis}}

The proof will compute all $S-$pairs of elements in the set $\Sigma=\Sigma(i_{1},\ldots,i_{j})$.
As usual, pairs $F,G$ such $\gcd({\rm in}(F),{\rm in}(G))=1$ will be overlooked.

\medskip

\noindent {\bf Case 1.} $S(K_{i,k},K_{i',k'}).$\\
In this case, $\mathrm{in}(K_{i,k})=x_{k}^{a}y_{i}$ and $\mathrm{in}(K_{i',k'})=x_{k'}^{a}y_{i'}.$

\smallskip

\noindent {\bf \small Case 1.1.} Let $i<k<i'<k'.$ No action here
since $\mathrm{in}(K_{i,k})$ and $\mathrm{in}(K_{i',k'})$ are relatively prime.

\smallskip

\noindent {\bf \small Case 1.2.} Let $i<i'<k<k'.$ No action here
since $\mathrm{in}(K_{i,k})$ and $\mathrm{in}(K_{i',k'})$ are relatively prime.

\smallskip

\noindent {\bf \small Case 1.3.} Let $i'<i<k<k'.$ No action here
since $\mathrm{in}(K_{i,k})$ and $\mathrm{in}(K_{i',k'})$ are relatively prime.

\smallskip

\noindent {\bf \small Case 1.4.} Let $k=k'$. Then $$S(K_{i,k},K_{i',k})= \frac{x_{k}^{a}y_{i}y_{i'}}{-x_{k}^{a}y_{i}}K_{i,k}-
\frac{x_{k}^{a}y_{i}y_{i'}}{-x_{k}^{a}y_{i'}}K_{i',k}
=-y_{k}(x_{i}^{a}y_{i'}-x_{i'}^{a}y_{i})\equiv 0\mod\Sigma.$$

\smallskip

\noindent {\bf \small Case 1.5.} Let $i<k=i'<k'.$  No action here
since $\mathrm{in}(K_{i,k})$ and $\mathrm{in}(K_{k,k'})$ are relatively prime.

\smallskip

\noindent {\bf \small Case 1.6.} Let $i'<k'=i<k.$ No action here
since $\mathrm{in}(K_{i,k})$ and $\mathrm{in}(K_{k,i})$ are relatively prime.

\smallskip

\noindent {\bf \small Case 1.7.} Let $i=i'.$ Then $$S(K_{i,k},K_{i,k'})= \frac{x_{k}^{a}x_{k'}^{a}y_{i}}{-x_{k}^{a}y_{i}}K_{i,k}-\frac{x_{k}^{a}x_{k'}^{a}y_{i}}{-x_{k'}^{a}y_{i}}K_{i,k'}=-x_{i}^{a}(x_{k'}^{a}y_{k}-x_{k}^{a}y_{k'})\equiv 0\mod\Sigma.$$

\smallskip

\noindent {\bf Case 2.} $S(L_{j},L_{j'}),$ with $j<j'.$\\ In this case, $\mathrm{in}(L_{j})=x_{j}^{a-b}w$ and $\mathrm{in}(L_{j'})=x_{j'}^{a-b}w.$\\ Then $$S(L_{j},L_{j'})=\frac{x_{j}^{a-b}x_{j'}^{a-b}w}{x_{j}^{a-b}w}L_{j}-\frac{x_{j}^{a-b}x_{j'}^{a-b}w}{x_{j'}^{a-b}w}L_{j'}=P(j,j')^bK_{j,j'}\equiv 0\mod\Sigma.$$

\smallskip

\noindent {\bf Case 3.} $S(L_{j},K_{i,k}).$\\
In this case, $\mathrm{in}(L_{j})=x_{j}^{a-b}w$ and $\mathrm{in}(K_{i,k})=x_{k}^{a}y_{i}.$

\smallskip

\noindent {\bf \small Case 3.1.} Let $j<i<k.$ No action here
since $\mathrm{in}(L_{j})$ and $\mathrm{in}(K_{i,k})$ are relatively prime.

\smallskip

\noindent {\bf \small Case 3.2.} Let $i<j<k.$ No action here
since $\mathrm{in}(L_{j})$ and $\mathrm{in}(K_{i,k})$ are relatively prime.

\smallskip

\noindent {\bf \small Case 3.3.} Let $i<k<j.$ No action here
since $\mathrm{in}(L_{j})$ and $\mathrm{in}(K_{i,k})$ are relatively prime.

\smallskip

\noindent {\bf \small Case 3.4.} Let $j=i<k.$ No action here
since $\mathrm{in}(L_{i})$ and $\mathrm{in}(K_{i,k})$ are relatively prime.

\smallskip

\noindent {\bf \small Case 3.5.} Let $i<j=k.$ Then $$S(L_{k},K_{i,k})=\frac{x_{k}^{a}wy_{i}}{x_{k}^{a-b}w}L_{k}-\frac{x_{k}^{a}wy_{i}}{-x_{k}^{a}y_{i}}K_{i,k}=x_{i}^{b}y_{k}L_{i}\equiv 0\mod\Sigma.$$

\smallskip

\noindent {\bf Case 4.} $S(K_{u,k},H_j^{i_1,\ldots,i_j})$. \\ In this case, $\mathrm{in}(K_{u,k})=x_{k}^{a}y_{u}$ and $\mathrm{in}(H_j^{i_1,\ldots,i_j})=(x_{i_{1}}\cdots x_{i_{j}})^{a-jb}w^{j}.$

\smallskip

\noindent {\bf \small Case 4.1.} Let $u<k$, $u,k \notin \{i_{1},\ldots,i_{j}\}$. No action here
since $\mathrm{in}(K_{u,k})$ and $\mathrm{in}(H_j^{i_1,\ldots,i_j})$ are relatively prime.

\smallskip

\noindent {\bf \small Case 4.2.} Let $u<k$, $u \in \{i_{1},\ldots,i_{j}\}$ and $k\notin \{i_{1},\ldots,i_{j}\}$. No action here
since $\mathrm{in}(K_{u,k})$ and $\mathrm{in}(H_j^{i_1,\ldots,i_j})$ are relatively prime.

\smallskip

\noindent {\bf \small Case 4.3.} Let $u<k$, $u \notin \{i_{1},\ldots,i_{j}\}$, and $k\in\{i_{1},\ldots,i_{j}\}$. Then $$S(K_{u,k},H_j^{i_1,\ldots,i_j})=\frac{x_{k}^{a}y_{u}(x_{i_{1}}\cdots \widehat{x_{k}}\cdots x_{i_{j}})^{a-jb}w^{j} }{-x_{k}^{a}y_{u}}K_{u,k}-\frac{x_{k}^{a}y_{u}(x_{i_{1}}\cdots \widehat{x_{k}}\cdots x_{i_{j}})^{a-jb}w^{j} }{(x_{i_{1}}\cdots x_{k}\cdots x_{i_{j}})^{a-jb}w^{j}}H_{j}^{i_{1},\ldots,i_{j}}$$ $$\hspace{1cm}=(-x_{u}^{jb}y_{k})H_{j}^{I'} \equiv 0\mod\Sigma,~ \mathrm{where}~ I'=(\{i_{1},\ldots,i_{j}\} \setminus \{k\})\cup \{u\}.$$

\smallskip

\noindent {\bf \small Case 4.4.} Let $u<k$, $u,k \in \{i_{1},\ldots,i_{j}\}$. Then $$S(K_{u,k},H_j^{i_1,\ldots,i_j})=\frac{x_{k}^{a}y_{u}(x_{i_{1}}\cdots \widehat{x_{k}}\cdots x_{i_{j}})^{a-jb}w^{j} }{-x_{k}^{a}y_{u}}K_{u,k}-\frac{x_{k}^{a}y_{u}(x_{i_{1}}\cdots \widehat{x_{k}}\cdots x_{i_{j}})^{a-jb}w^{j} }{(x_{i_{1}}\cdots x_{k}\cdots x_{i_{j}})^{a-jb}w^{j}}H_{j}^{i_{1},\ldots,i_{j}}$$ $$\hspace{2,8cm}=-y_{k}[x_{u}^{a}(x_{i_{1}}\cdots \widehat{x_{k}}\cdots x_{i_{j}})^{a-jb}w^{j}-x_{k}^{jb}y_{u}P(i_{1},\ldots,i_{j})^{jb}y_{i_{1}}\cdots y_{u}\cdots \widehat{y_{k}} \cdots y_{i_{j}}]. $$
Since $x_u^{a-b}w=L_u+P(u)^by_u$, it obtains
$$S(K_{u,k},H_j^{i_1,\ldots,i_j})\equiv -P(i_{1},\ldots,\widehat{k},\ldots,i_{j})^{b} y_{u}y_{k} H_{j-1}^{I''}\equiv 0\mod\Sigma,$$ $\mathrm{where}~ I''=\{i_{1},\ldots,i_{j}\} \setminus \{k\}.$\\

\smallskip

\noindent {\bf Case 5.} $S(L_{u},H_j^{i_1,\ldots,i_j})$. \\
In this case, $\mathrm{in}(L_{u})=x_{u}^{a-b}w$ and $\mathrm{in}(H_j^{i_1,\ldots,i_j})=(x_{i_{1}}\cdots x_{i_{j}})^{a-jb}w^{j}.$

\smallskip

\noindent {\bf \small Case 5.1.} Let $u \notin \{i_{1},\ldots,i_{j}\}$. Then \begin{eqnarray}\nonumber
S(L_{u},H_j^{i_1,\ldots,i_j}) &=& -P(i_{1},\ldots,i_{j},u)^{b}[(x_{i_{1}}\cdots x_{i_{j}})^{a-(j-1)b}w^{j-1}y_{u}\\ \nonumber
&-& x_{u}^{a+(j-1)b}P(i_{1},\ldots,i_{j},u)^{(j-1)b}y_{i_{1}},\ldots,y_{i_{j}}].
\end{eqnarray}
Pick any subset $I'\subset \{i_{1},\ldots,i_{j}\}$, with $|I'|=j-1$ and reduce modulo $H_{j-1}^{I'}$ the monomial with $w$ occurring within the square brackets. The result is a binomial not involving $w$. By the same argument as before, we conclude that this pair reduces to $0$ modulo $\Sigma$.

\smallskip

\noindent {\bf \small Case 5.2.} Let $u \in \{i_{1},\ldots,i_{j}\}$. Then $$S(L_{u},H_j^{i_1,\ldots,i_j})=-P(i_{1},\ldots,i_{j})^{b}y_{u}H_{j-1}^{i_{1},\ldots,\widehat{u},\ldots,i_{j}}\equiv 0\mod\Sigma .$$

\smallskip

\noindent {\bf Case 6.} Consider $H_m^{i_1,\ldots,i_m}$ and $H_{m'}^{j_1,\ldots,j_{m'}}$, with the respective external degrees $m\leq m'$. Denote $\mathfrak{I}:=\{i_1,\ldots,i_m\},~\mathfrak{I}':=\{j_1,\ldots,j_{m'}\}$ and let $\mathfrak{I}\cap \mathfrak{I}'=\{k_{1},\ldots,k_{s}\}$, for some $s\in \{0,\ldots, m\}.$

Under the given order, the two leading terms of the two binomials are \\${\rm in}(H_m^{i_{1},\ldots,i_{m}})=(x_{i_{1}}\cdots x_{i_{m}})^{a- mb}w^m$ and ${\rm in}(H_{m'}^{j_{1},\ldots,j_{m'}})=(x_{j_{1}}\cdots x_{j_{m'}})^{a-m'b}w^{m'}$, so their least common multiple is $w^{m'}(x_{j_{1}}\cdots \widehat{x_{i_{1}}}\cdots \widehat{x_{i_{m}}} \cdots x_{j_{m'}})^{a- m'b}(x_{i_{1}}\cdots x_{i_{m}} )^{a-mb}$. Therefore
\begin{eqnarray}\nonumber
\nonumber
\lefteqn{S(H_m^{i_{1},\ldots,i_{m}},H_{m'}^{j_{1},\ldots,j_{m'} })=-P(i_{1},\ldots,j_{1},\ldots,i_{m},\ldots,j_{m'})^{mb}y_{k_{1}}\cdots y_{k_{s}} }\kern2cm\\ \nonumber
&&\cdot \left[w^{m'-m}(x_{j_{1}} \cdots\widehat{x_{i_{1}}} \cdots \widehat{x_{i_{m}}}\cdots x_{j_{m'}})^{a-m'b+mb}y_{i_{1}}\cdots\widehat{y_{j_{1}}}\cdots \widehat{y_{j_{m'}}} \cdots y_{i_m}\right.\\ \nonumber
&&-(x_{i_{1}}\cdots \widehat{x_{j_{1}}}\cdots \widehat{x_{j_{m'}}} \cdots x_{i_{m}} )^{a-mb+m'b}(x_{k_{1}} \cdots x_{k_{s}})^{(m'-m)b}\\ \nonumber
&&\cdot \left. P(i_{1},\ldots,j_{1},\ldots,i_{m},\ldots,j_{m'})^{(m'-m)b}y_{j_{1}}\cdots\widehat{y_{i_{1}}}\cdots \widehat{y_{i_{m}}} \cdots y_{j_{m'}} \right].
\end{eqnarray}
If $m'=m$, then the binomial inside the square brackets does not involve $w$ and therefore reduces to $0$ modulo $\Sigma$ by previous cases.

If $m'>m$, then $|\mathfrak{I}'|=m'>m=|\mathfrak{I}|$, and therefore $|\mathfrak{I}'\setminus \mathfrak{I}|\geq m'-m\geq 1$. Let $\widetilde{\mathfrak{I}}\subset \mathfrak{I}'\setminus \mathfrak{I}$ with $|\widetilde{\mathfrak{I}}|=m'-m$. Reducing the binomial inside the square brackets modulo $H_{m'-m}^{\widetilde{\mathfrak{I}}}\in\Sigma$, will result in the cancellation of the monomial involving $w$, hence we are back to the previous situation.
\qed

\subsection{Proof of Proposition~\ref{colon-of-initials}}
To apply Lemma~\ref{lemma_colon} (i) we set ourselves to compute a minimal set of generators of the intersection of the two initial ideals on the left hand side, then divide each generator by the initial term of $H_{j'}^{k_{1},\ldots,k_{j'}}$.
A minimal set of generators of the intersection turns out to be the set of the least common multiples of ${\rm in}(H_{j'}^{k_{1},\ldots,k_{j'}})$ and each minimal generator of ${\rm in}({\mathcal H}(i_{1},\ldots,i_{j}))$.

We separate the two cases, according as to whether $j'=j$ or $j'=j+1$.

\medskip

(a) {\sc Same degree:} $j=j'$

One has ${\rm in}(H_{j}^{k_{1},\ldots,k_{j}})=(x_{k_{1}}\cdots x_{k_{j}})^{a-jb}w^{j}$.
Drawing upon (\ref{initial_ideal}), according to the external degree of a monomial, we have

\smallskip

\textbf{Degree 1:}
\begin{itemize}

\item[\textbullet] $x_{d}^{a-b}w,~d\notin\{k_{1},\ldots,k_{j}\}$ and $d<k_{j}$ (coming from $L_{d}\in \mathcal{L}$)
$$\frac{\mathrm{lcm}(x_{d}^{a-b}w,(x_{k_{1}}\cdots x_{k_{j}})^{a-jb}w^{j})}{(x_{k_{1}}\cdots x_{k_{j}})^{a-jb}w^{j}}=x_{d}^{a-b}.$$
As $j\leq p-1$, and $a>(p-1)b$, then $a-jb>0$. But then $x_{d}^{a-b}=x_{d}^{(j-1)b}x_{d}^{a-jb}$, and $x_{d}^{a-jb}$ is among the generators listed in the right hand side monomial ideal.

\item[\textbullet] $x_{d}^{a-b}w,~d\notin\{k_{1},\ldots,k_{j}\}$ and $d>k_{j}$ (coming from $L_{d}\in \mathcal{L}$)
$$\frac{\mathrm{lcm}(x_{d}^{a-b}w,(x_{k_{1}}\cdots x_{k_{j}})^{a-jb}w^{j})}{(x_{k_{1}}\cdots x_{k_{j}})^{a-jb}w^{j}}=x_{d}^{a-b}.$$

As above, $x_{d}^{a-jb}$ is among the generators listed in the right hand side monomial ideal.

\item[\textbullet] $x_{d}^{a-b}w,~d\in\{k_{1},\ldots,k_{j}\}$ (coming from $L_{d}\in \mathcal{L}$)

$$\frac{\mathrm{lcm}(x_{d}^{a-b}w,(x_{k_{1}}\cdots x_{k_{j}})^{a-jb}w^{j})}{(x_{k_{1}}\cdots x_{k_{j}})^{a-jb}w^{j}}=x_{d}^{(j-1)b},$$
which is among the generators listed in the right hand side monomial ideal.

\item[\textbullet] $x_{d}^{a}y_{v},~d\notin\{k_{1},\ldots,k_{j}\}$ and $d<k_{j}$ (coming from $K_{d,v}\in \mathcal{L}$)
$$\frac{\mathrm{lcm}(x_{d}^{a}y_{v},(x_{k_{1}}\cdots x_{k_{j}})^{a-jb}w^{j})}{(x_{k_{1}}\cdots x_{k_{j}})^{a-jb}w^{j}}=x_{d}^{a}y_{v}$$
Note that $x_{d}^{a}y_{v}=(x_{d}^{jb}y_{v})x_{d}^{a-jb}$, while $x_{d}^{a-jb}$ is among the generators listed in the right hand side monomial ideal.

\item[\textbullet] $x_{d}^{a}y_{v},~d\notin\{k_{1},\ldots,k_{j}\}$ and $d>k_{j}$ (coming from $K_{d,v}\in \mathcal{L}$)
$$\frac{\mathrm{lcm}(x_{d}^{a}y_{v},(x_{k_{1}}\cdots x_{k_{j}})^{a-jb}w^{j})}{(x_{k_{1}}\cdots x_{k_{j}})^{a-jb}w^{j}}=x_{d}^{a}y_{v}$$
One has $x_{d}^{a}y_{v}=(x_{d}^{b}y_{v})x_{d}^{a-b}$, while $x_{d}^{a-b}$ is among the generators listed in the right hand side monomial ideal.

\item[\textbullet] $x_{d}^{a}y_{v},~d\in\{k_{1},\ldots,k_{j}\}$ (coming from $K_{d,v}\in \mathcal{L}$)
$$\frac{\mathrm{lcm}(x_{d}^{a}y_{v},(x_{k_{1}}\cdots x_{k_{j}})^{a-jb}w^{j})}{(x_{k_{1}}\cdots x_{k_{j}})^{a-jb}w^{j}} =x_{d}^{jb}y_{v}$$
Note that $x_{d}^{jb}y_{v}=(x_{d}^{(j-1)b}y_{v})x_{d}^{b}$,  so once more we get a generator listed in the right hand side monomial ideal.

\end{itemize}

\medskip

\textbf{Degree} {${\mathbf s}\,  (2\leq s \leq j-1)$:}

\begin{itemize}
\item[\textbullet] $(x_{d_{1}}\cdots x_{d_{r}}x_{q_{1}}\cdots  x_{q_{s-r}})^{a-sb}w^{s},~\{d_{1}<\ldots <d_{r}\}\cap\{k_{1},\ldots,k_{j}\}=\emptyset$, $d_1<k_{j}$ and $\{q_{1}<\cdots<q_{s-r}\}\subset\{k_{1},\ldots,k_{j}\}$.
\begin{eqnarray}\nonumber
&&\kern-0.6cm \frac{\mathrm{lcm}((x_{d_{1}}\cdots x_{d_{r}}x_{q_{1}}\cdots  x_{q_{s-r}})^{a-sb}w^{s},(x_{k_{1}}\cdots x_{k_{j}})^{a-jb}w^{j})} {(x_{k_{1}}\cdots x_{k_{j}})^{a-jb}w^{j}}\\ \nonumber
&=&
(x_{d_{1}}\cdots x_{d_{r}})^{a-sb} (x_{q_{1}}\cdots  x_{q_{s-r}})^{(j-s)b}
\end{eqnarray}
Note that $x_{d_1}^{a-sb}$ is a factor thereof factoring further as
 $x_{d_1}^{a-sb}=x_{d_1}^{(j-s)b}x_{d_1}^{a-jb}$, while $x_{d_1}^{a-jb}$ is among the generators listed in the right hand side monomial ideal since $d_1\notin \{k_1,\ldots,k_j\}$ and  $d_1<k_j$.

\item[\textbullet] $(x_{q_{1}}\cdots  x_{q_{s}})^{a-sb}w^{s},~\{q_{1}<\cdots<q_{s}\}\subset\{k_{1},\ldots,k_{j}\}$.

$$\frac{\mathrm{lcm}((x_{q_{1}}\cdots  x_{q_{s}})^{a-sb}w^{s},(x_{k_{1}}\cdots x_{k_{j}})^{a-jb}w^{j})}{(x_{k_{1}}\cdots x_{k_{j}})^{a-jb}w^{j}}=(x_{q_{1}}\cdots  x_{q_{s}})^{(j-s)b}$$

\item[\textbullet] $(x_{q_{1}}\cdots  x_{q_{r}}x_{d_{1}}\cdots x_{d_{s-r}})^{a-sb}w^{s},$ $\{q_{1}<\cdots<q_{r}\}\subset\{k_{1},\ldots,k_{j}\}, ~d_{1}<\ldots<d_{s-r}$ with $k_{j}<d_{1}$.
\begin{eqnarray}\nonumber
&&\kern-0.6cm \frac{\mathrm{lcm}((x_{q_{1}}\cdots  x_{q_{r}}x_{d_{1}}\cdots x_{d_{s-r}})^{a-sb}w^{s},(x_{k_{1}}\cdots x_{k_{j}})^{a-jb}w^{j})}
{(x_{k_{1}}\cdots x_{k_{j}})^{a-jb}w^{j}}\\ \nonumber
&=&(x_{q_{1}}\cdots  x_{q_{r}})^{(j-s)b}(x_{d_{1}}\cdots x_{d_{s-r}})^{a-sb}
\end{eqnarray}

\item[\textbullet] $(x_{d_{1}}\cdots x_{d_{s}})^{a-sb}w^{s},~d_{1}<\ldots<d_{s}$ with $k_{j}<d_{1}$.
$$\frac{\mathrm{lcm}((x_{d_{1}}\cdots x_{d_{s}})^{a-sb}w^{s},(x_{k_{1}}\cdots x_{k_{j}})^{a-jb}w^{j})}{(x_{k_{1}}\cdots x_{k_{j}})^{a-jb}w^{j}}=(x_{d_{1}}\cdots x_{d_{s}})^{a-sb}.$$
In all three cases above the resulting monomial is among the generators listed in the right hand side monomial ideal.

\end{itemize}

\textbf{Degree} $\bf j:$
\begin{itemize}
\item[\textbullet] $(x_{k_{1}}\cdots \widehat{x_{k_{c}}}\ldots x_{k_{j}})^{a-jb}{x_d}^{a-jb}w^{j},d \notin\{k_{1},\ldots,k_{j}\}$, $d<k_{j}$, $c\in\{1,\ldots,j\}$.

$$\frac{\mathrm{lcm}((x_{k_{1}}\cdots \widehat{x_{k_{c}}}\ldots x_{k_{j}})^{a-jb}{x_d}^{a-jb}w^{j},(x_{k_{1}}\cdots x_{k_{j}})^{a-jb}w^{j})}{(x_{k_{1}}\cdots x_{k_{j}})^{a-jb}w^{j}}=x_{d}^{a-jb}.$$
Again we conclude as before.

\end{itemize}

\medskip

(b) {\sc Degree jump:} $j'=j+1$

\medskip

We now consider the case where the degree goes up, that is, one is dealing with  $${\rm in}(H_{j+1}^{1,\ldots,j+1})=(x_1\cdots x_{j+1})^{a-(j+1)b} w^{j+1}.$$

\medskip

We go through similar calculations as before. In each case below  the resulting monomial is among the generators listed in the right hand side monomial ideal.

\bigskip

\textbf{Degree 1:}
\begin{itemize}
\item[\textbullet] $x_{d}^{a-b}w,~d\notin\{1,\ldots,j+1\}$.
$$\frac{\mathrm{lcm}(x_{d}^{a-b}w,(x_{1}\cdots x_{j+1})^{a-(j+1)b}w^{j+1})}{(x_{1}\cdots x_{j+1})^{a-(j+1)b}w^{j+1}}=x_{d}^{a-b}$$

\item[\textbullet] $x_{d}^{a-b}w,~d\in\{1,\ldots,j+1\}$.
$$\frac{\mathrm{lcm}(x_{d}^{a-b}w,(x_{1}\cdots x_{j+1})^{a-(j+1)b}w^{j+1})}{(x_{1}\cdots x_{j+1})^{a-(j+1)b}w^{j+1}}=x_{d}^{jb}$$

\item[\textbullet] $x_{d}^{a}y_{v},~d\notin\{1,\ldots,j+1\}$.
$$\frac{\mathrm{lcm}(x_{d}^{a}y_{v},(x_{1}\cdots x_{j+1})^{a-(j+1)b}w^{j+1})}{(x_{1}\cdots x_{j+1})^{a-(j+1)b}w^{j+1}}=x_{d}^{a}y_{v}$$
Note that $x_{d}^{a}y_{v}=(x_{d}^{b}y_{v})x_{d}^{a-b}$.

\item[\textbullet] $x_{d}^{a}y_{v},~d\in\{1,\ldots,j+1\}$.
$$\frac{\mathrm{lcm}(x_{d}^{a}y_{v},(x_{1}\cdots x_{j+1})^{a-(j+1)b}w^{j+1})}{(x_{1}\cdots x_{j+1})^{a-(j+1)b}w^{j+1}}=x_{d}^{(j+1)b}y_{v}$$
Note that $x_{d}^{(j+1)b}y_{v}=(x_{d}^{b}y_{v})x_{d}^{jb}$.

\end{itemize}

\textbf{Degree} {${\mathbf s}\,  (2\leq s \leq j)$:}

\begin{itemize}

\item[\textbullet] $(x_{q_{1}}\cdots  x_{q_{s}})^{a-sb}w^{s},$ with $\{q_{1}<\cdots<q_{s}\}\subset\{1,\ldots,j+1\}$.
$$\frac{\mathrm{lcm}((q_{1}\cdots  x_{q_{s}})^{a-sb}w^{s},(x_{1}\cdots x_{j+1})^{a-(j+1)b}w^{j+1})}{(x_{1}\cdots x_{j+1})^{a-(j+1)b}w^{j+1}}=(x_{q_{1}}\cdots  x_{q_{s}})^{(j+1-s)b}$$

\item[\textbullet] $(x_{q_{1}}\cdots  x_{q_{r}}x_{d_{1}}\cdots x_{d_{s-r}})^{a-sb}w^{s},$ $\{q_{1}<\cdots<q_{r}\}\subset\{1,\ldots,j+1\}, ~d_{1}<\ldots<d_{s-r}$ with $j+1<d_{1}$.
\begin{eqnarray}\nonumber
&&\kern-0.6cm\frac{\mathrm{lcm}((x_{q_{1}}\cdots  x_{q_{r}}x_{d_{1}}\cdots x_{d_{s-r}})^{a-sb}w^{s},(x_{1}\cdots x_{j+1})^{a-(j+1)b}w^{j+1})}{(x_{1}\cdots x_{j+1})^{a-(j+1)b}w^{j+1}}\\ \nonumber
&=&
(x_{q_{1}}\cdots  x_{q_{r}})^{(j+1-s)b}(x_{d_{1}}\cdots x_{d_{s-r}})^{a-sb}
\end{eqnarray}

\item[\textbullet] $(x_{d_{1}}\cdots x_{d_{s}})^{a-sb}w^{s},~d_{1}<\ldots<d_{s-r}$ with $j+1<d_{1}$.
$$\frac{\mathrm{lcm}((x_{d_{1}}\cdots x_{d_{s}})^{a-sb}w^{s},(x_{1}\cdots x_{j+1})^{a-(j+1)b}w^{j+1})}{(x_{1}\cdots x_{j+1})^{a-(j+1)b}w^{j+1}}=(x_{d_{1}}\cdots x_{d_{s}})^{a-sb}$$

\end{itemize}

To conclude the present case of degree jump, we stress the limit situation where the degree jumps to the highest possible degree of a Sylvester form.
It is convenient to separate the two basic settings:

\medskip

{\sc Setting $a>nb$.}

The expected outcome is  ${\rm in}(\mathcal{H}(2,\ldots,n):{\rm in}(H_{n}^{1,\ldots,n})=(x_{1},\ldots,x_{n})^{(n-1)b}S$
and the calculation of the required least common multiples is included in the general calculation above, setting $j=n$.

\smallskip

{\sc Setting $a\leq nb,~(p-1)b<a\leq pb$.}

Here
$H_{p}^{1,\ldots,p}=w^{p}-(x_{p+1}\cdots x_{n})^{pb}(x_{1}\cdots x_{p})^{pb-a}y_{1}\cdots y_{p}$.
The typical expected generator has one of the following forms
 $$(x_{r_{1}}\cdots x_{r_{s}})^{a-sb} \; {\rm and}\;\; (x_{q_1}\cdots x_{q_{r}}x_{d_{1}}\cdots x_{d_{s-r}} )^{a-sb},
 $$
where $s\in \{2,\ldots,p-1\},~r\in \{0,\ldots,s-1\},~\{r_{1}<\cdots<r_{s}\}\subset\{1,\ldots,p\},~\{q_{1}<\cdots<q_{r}\}\subset\{1,\ldots,p\} ,~\{d_{1}<\cdots<d_{s-r}\}\subset\{p+1,\ldots,n\}$.

\smallskip

Here is the calculation for this setting, according to the external degrees of the generating monomials:

\smallskip

\textbf{Degree 1:}

\begin{itemize}
\item[\textbullet] $x_{d}^{a-b}w,~d\notin\{1,\ldots,p\}$.
$$\frac{\mathrm{lcm}(x_{d}^{a-b}w,w^{p})}{w^{p}}=x_{d}^{a-b}$$

\item[\textbullet] $x_{d}^{a-b}w,~d\in\{1,\ldots,p\}$.
$$\frac{\mathrm{lcm}(x_{d}^{a-b}w,w^{p})}{w^{p}}=x_{d}^{a-b}$$

\item[\textbullet] $x_{d}^{a}y_{v},~d\notin\{1,\ldots,p\}$.
$$\frac{\mathrm{lcm}(x_{d}^{a}y_{v},w^{p})}{w^{p}}=x_{d}^{a}y_{v}$$
Note that $x_{d}^{a}y_{v}=(x_{d}^{b}y_{v})x_{d}^{a-b}$.

\item[\textbullet] $x_{d}^{a}y_{v},~d\in\{1,\ldots,p\}$.
$$\frac{\mathrm{lcm}(x_{d}^{a}y_{v},w^{p})}{w^{p}}=x_{d}^{a}y_{v}$$
Note that $x_{d}^{a}y_{v}=(x_{d}^{b}y_{v})x_{d}^{a-b}$.

\end{itemize}

\smallskip

\textbf{Degree} ${\mathbf s}\, (2\leq {\mathbf s} \leq p-1$):
\begin{itemize}
\item[\textbullet] $(x_{q_{1}}\cdots  x_{q_{s}})^{a-sb}w^{s},~ \{q_{1}<\cdots<q_{s}\}\subset\{1,\ldots,p\}$ \\
$$\frac{\mathrm{lcm}((q_{1}\cdots  x_{q_{s}})^{a-sb}w^{s},w^{p})}{w^{p}}=(x_{q_{1}}\cdots  x_{q_{s}})^{a-sb}$$

\item[\textbullet] $(x_{q_{1}}\cdots  x_{q_{r}}x_{d_{1}}\cdots x_{d_{s-r}})^{a-sb}w^{s},~\{q_{1}<\cdots<q_{s}\}\subset\{1,\ldots,p\}, ~d_{1}<\ldots<d_{s-r}$ with $p<d_{1}$  \\
$$\frac{\mathrm{lcm}((x_{q_{1}}\cdots  x_{q_{r}}x_{d_{1}}\cdots x_{d_{s-r}})^{a-sb}w^{s},w^{p})}{w^{p}}=(x_{q_{1}}\cdots  x_{q_{r}}x_{d_{1}}\cdots x_{d_{s-r}})^{a-sb}$$

\item[\textbullet] $(x_{d_{1}}\cdots x_{d_{s}})^{a-sb}w^{s},~d_{1}<\ldots<d_{s}$ with $p<d_{1}$  \\
$$\frac{\mathrm{lcm}((x_{d_{1}}\cdots x_{d_{s}})^{a-sb}w^{s},w^{p})}{w^{p}}=(x_{d_{1}}\cdots x_{d_{s}})^{a-sb}$$

\end{itemize}

\smallskip

This concludes the proof  of the proposition.

\subsection{Proof of Proposition~\ref{successive_colon}}
We just have to prove the inclusion $\supset$ since the inclusion $\subset$ follows from it by applying Lemma~\ref{lemma_colon} (ii).

Again, we deal with two cases, according to the established sets of generators for the right hand side of the stated equality
in either (a) or (b) of  Proposition~\ref{G-basis}.

\subsubsection{Same degree}

\begin{itemize}
\item[\textbullet] $x_{s}^{a-b},~s=k_{j}+1,\ldots,n.$

$$x_{s}^{a-b}H_{j}^{k_{1},\ldots,k_{j}}=(x_{k_{1}}\cdots x_{k_{j}})^{a-jb}w^{j-1}L_{s} + x_{k_{1}}^{a-(j-1)b}P(k_1,\ldots,k_j,s)^{b}y_{s}H_{j-1}^{k_{2},\ldots,k_{j}} $$ $$+ P(k_1,\ldots,k_j,s)^{jb}x_{s}^{(j-1)b}y_{k_{2}}\cdots y_{k_{j}}K_{k_{1},s}.$$

\item[\textbullet] $x_{k_{s}}^{(j-1)b},~s=1,\ldots,j.$
$$x_{k_{s}}^{(j-1)b}H_{j}^{k_{1},\ldots,k_{j}}=(x_{k_1} \cdots \widehat{x_{k_s}} \cdots x_{k_j})^{a-jb}w^{j-1}L_{k_{s}}\\ \nonumber
+P(k_1,\ldots,k_s,\ldots, k_j)^{b}y_{k_{s}}H_{j-1}^{k_{1},\ldots,\widehat{k_{s}},\ldots,k_{j} }.
$$

\item[\textbullet] $(x_{r_{1}}\cdots x_{r_{s}})^{(j-s)b},~s\in\{1,\ldots,j-1\},~\{r_{1}<\cdots<r_{s}\}\subset \{k_{1},\ldots,k_{j}\} .$

$$(x_{r_{1}}\cdots x_{r_{s}})^{(j-s)b}H_{j}^{k_{1},\ldots,k_{j}}= (x_{k_1} \cdots \widehat{x_{r_1}} \cdots \widehat{x_{r_s}} \cdots x_{k_j})^{a-jb}w^{j-s}H_{s}^{r_{1},\ldots,r_{s}}+$$ $$P(k_{1},\ldots, r_1,\ldots,r_s,\ldots,k_{j})^{sb}y_{r_{1}}\cdots y_{r_{s}}H_{j-s}^{k_{1},\ldots,\widehat{r_{1}},\ldots,\widehat{r_{s}},\ldots,k_{j}}  $$

\item[\textbullet] $x_{r}^{a-jb},~r<k_{1}$.

$$x_{r}^{a-jb}H_{j}^{k_{1},\ldots,k_{j}}=x_{k_{1}}^{a-jb}H_{j}^{r,k_{2},\ldots,k_{j}}-P(r,k_1,\ldots,k_j)^{jb} y_{k_{2}}\cdots y_{k_{j}}K_{r,k_{1}}. $$

\item[\textbullet] $x_{r}^{a-jb},k_{i}<r<k_{i+1},~i=1,\ldots,j-1$~.

\begin{eqnarray}\nonumber
x_{r}^{a-jb}H_{j}^{k_{1},\ldots,k_{j}}&=&x_{k_{i+1}}^{a-jb}H_{j}^{k_{1},k_{2},\ldots,k_i,r,\widehat{k_{i+1}},\ldots, k_{j}}\\ \nonumber
&-&P(k_1,\ldots,k_i,r,k_{i+1},\ldots,k_j)^{jb} y_{k_{1}}\cdots y_{k_{i}}\widehat{y_{k_{i+1}}} y_{k_{i+2}}\cdots y_{k_{j}} K_{r,k_{i+1}}.
\end{eqnarray}

\item[\textbullet] $(x_{q_{1}}\cdots x_{q_{r}})^{(j-s)b}(x_{d_{1}}\cdots x_{d_{s-r}})^{a-sb}$.
\end{itemize}

Since this case is a lot more involved than the previous ones, we chose to formulate it as a lemma.

\begin{Lemma}\label{last_case}
Fix an integer $2\leq j\leq p-1$ and an ordered subset  $\{k_1,\ldots,k_j\}\subset \{1,\ldots,n\}$.
Let there be given integers $s\in \{1,\ldots,j-1\},\, r\in \{0,\ldots,s-1\}$
and ordered subsets $\{q_1,\ldots ,q_r\}\subset \{k_1,\ldots,k_j\}$ and $\{d_1,\ldots ,d_{s-r}\}\subset \{1,\ldots,n\}\setminus \{k_1,\ldots,k_j\}$,  with $k_j<d_1$.
Consider a $2$-partition of $\{k_1,\ldots,k_j\}\setminus \{q_1,\ldots ,q_r\}$  by ordered subsets $\{k_{m_{1}},\ldots,k_{m_{j-s}}\}$ and $\{n_{1},\ldots,n_{s-r}\}$.
Set
$$Q:=P(k_1,\ldots ,k_j,\, d_1,\ldots,d_{s-r})^{jb}\,(x_{q_{1}}\cdots x_{q_{r}})^{(j-s)b}y_{q_{1}}\cdots y_{q_{r}}y_{k_{m_{1}}}\cdots y_{k_{m_{j-s}}}.$$
Then
\begin{eqnarray}\nonumber
\lefteqn{ \kern-1.3cm  (x_{q_{1}}\cdots x_{q_{r}})^{(j-s)b}(x_{d_{1}}\cdots x_{d_{s-r}})^{a-sb}H_{j}^{k_{1},\ldots,k_{j}}=
(x_{k_{m_{1}}}\cdots x_{k_{m_{j-s}}}\, x_{n_{1}}\cdots x_{n_{s-r}})^{a-jb}  } \kern-0.5cm\\ \nonumber
&&\cdot w^{j-s}H_{s}^{q_{1},\ldots,q_{r},d_{1},\ldots,d_{s-r}}\\ \nonumber
&&+P(k_1,\ldots,k_j,d_1,\ldots ,d_{s-r})^{sb}
\cdot (x_{n_{1}}\cdots x_{n_{s-r}})^{a-(j-s)b} \\\nonumber
&&\cdot y_{q_{1}}\cdots y_{q_{r}}y_{d_{1}}\cdots y_{d_{s-r}}H_{j-s}^{k_{m_{1}},\ldots,k_{m_{j-s}}}\\\nonumber
&&+ \sum_{c=1}^{s-r}\, (x_{n_{c+1}}\cdots x_{n_{s-r}} )^{a}(x_{d_{1}}\cdots x_{d_{c-1}} )^{a+(j-s)b}(x_{d_{c}}\cdots x_{d_{s-r}} )^{(j-s)b}\,Q\\\nonumber
&&\cdot y_{d_{c+1}}\cdots y_{d_{s-r}} y_{n_{1}}\cdots y_{n_{c-1}}K_{n_{c},d_{c}},
\end{eqnarray}
with the convention that $x_{d_0}=y_{n_0}=x_{n_{s-r+1}}=x_{d_{s-r+1}}=1$.
\end{Lemma}
\demo
Although the above expression is verifiable by expanding the right hand side, the idea to get at it is by no means obvious. Since similar expressions will appear in the sequel, we will now explain its main core.
Thus, first write
\begin{eqnarray}\nonumber
\lefteqn{ 
 (x_{q_{1}}\cdots x_{q_{r}})^{(j-s)b}(x_{d_{1}}\cdots x_{d_{s-r}})^{a-sb}H_{j}^{k_{1},\ldots,k_{j}}=
(x_{q_{1}}\cdots x_{q_{r}})^{(j-s)b}(x_{d_{1}}\cdots x_{d_{s-r}})^{a-sb}  } \kern1.2cm\\ \nonumber
&&\cdot\left( (x_{k_1}\cdots x_{k_j})^{a-jb}w^j-P(k_1,\ldots,k_j)^{jb}\,y_{k_1}\cdots y_{k_j}  \right)\\ 
&&=(x_{k_1}\cdots \widehat{x_{q_{1}}}\cdots \widehat{x_{q_{r}}}\cdots x_{k_j})^{a-jb}\,(x_{q_{1}}\cdots x_{q_{r}})^{a-sb}\,(x_{d_{1}}\cdots x_{d_{s-r}})^{a-sb}w^j\\ 
&&-(x_{q_{1}}\cdots x_{q_{r}})^{(j-s)b}\,(x_{d_{1}}\cdots x_{d_{s-r}})^{a+(j-s)b}\, P(k_1,\ldots,k_j,d_1,\ldots,d_{s-r})^{jb}y_{k_1}\cdots y_{k_j}
\end{eqnarray}

Next, we rewrite each of the numbered expressions above.
Using the partition explained above, one can write

\begin{eqnarray}\nonumber
\lefteqn{(x_{k_{m_{1}}}\cdots x_{k_{m_{j-s}}}\, x_{n_{1}}\cdots x_{n_{s-r}})^{a-jb}\,w^{j-s} H_s^{q_1,\ldots,q_r,\, d_1,\ldots,d_{s-r}} }       \kern1cm\\ \nonumber
&&=(x_{k_{m_{1}}}\cdots x_{k_{m_{j-s}}}\, x_{n_{1}}\cdots x_{n_{s-r}})^{a-jb}\,w^{j-s}\left((x_{q_{1}}\cdots x_{q_{r}}\,x_{d_{1}}\cdots x_{d_{s-r}})^{a-sb}w^s\right.\\ \nonumber
&&\left. -P(q_1,\ldots,q_r,d_1,\ldots,d_{s-r})^{sb} y_{q_1}\cdots y_{q_r}\,y_{d_1}\cdots y_{d_{s-r}}\right)\\ 
&&=(x_{k_1}\cdots \widehat{x_{q_{1}}}\cdots \widehat{x_{q_{r}}}\cdots x_{k_j})^{a-jb}\,(x_{q_{1}}\cdots x_{q_{r}})^{a-sb}\,(x_{d_{1}}\cdots x_{d_{s-r}})^{a-sb}w^j \\ 
&&- P(k_1,\ldots,k_j, d_1,\ldots,d_{s-r})^{sb}\,(x_{k_{m_{1}}}\cdots x_{k_{m_{j-s}}}\, x_{n_{1}}\cdots x_{n_{s-r}})^{a-(j-s)b}\,w^{j-s}\\ \nonumber
&&\cdot\; y_{q_1}\cdots y_{q_r}\,y_{d_1}\cdots y_{d_{s-r}}.
\end{eqnarray}

The first numbered expression above is exactly the first numbered expression in the previous display.
However, the second numbered expression above does not coincide with the second numbered expression in the previous display, so there is a little more to pursue in order to cancel this expression by bringing up an expression involving another Sylvester form:

\begin{eqnarray}\nonumber
\lefteqn{
P(k_1,\ldots,k_j, d_1,\ldots,d_{s-r})^{sb}\,(x_{n_{1}}\cdots x_{n_{s-r}})^{a-(j-s)b} y_{q_1}\cdots y_{q_r}\,y_{d_1}\cdots y_{d_{s-r}}H_{j-s}^{k_{m_1}\ldots, k_{m_{j-s}}} } \kern1cm\\ \nonumber
&&= P(k_1,\ldots,k_j, d_1,\ldots,d_{s-r})^{sb}\,(x_{n_{1}}\cdots x_{n_{s-r}})^{a-(j-s)b} y_{q_1}\cdots y_{q_r}\,y_{d_1}\cdots y_{d_{s-r}}\\ \nonumber
&&\cdot\, \left(  (x_{k_{m_1}}\cdots x_{k_{m_{j-s}}})^{a-(j-s)b}\,w^{j-s}
- P(k_{m_1},\ldots, k_{m_{j-s}})^{(j-s)b}\,y_{k_{m_{1}}}\cdots y_{k_{m_{j-s}}}\right)\\ \nonumber
&&= P(k_1,\ldots,k_j, d_1,\ldots,d_{s-r})^{sb}\,(x_{k_{m_{1}}}\cdots x_{k_{m_{j-s}}}\, x_{n_{1}}\cdots x_{n_{s-r}})^{a-(j-s)b}\,w^{j-s}\\ \nonumber
&&\cdot\; y_{q_1}\cdots y_{q_r}\,y_{d_1}\cdots y_{d_{s-r}}\\ \nonumber
&&-  P(k_1,\ldots,k_j, d_1,\ldots,d_{s-r})^{jb}\, (x_{q_1}\cdots x_{q_r})^{(j-s)b} \, (x_{n_{1}}\cdots x_{n_{s-r}})^{a}\, (x_{d_1}\cdots x_{d_{s-r}})^{(j-s)b}\\ \nonumber
&&\cdot\, y_{q_1}\cdots y_{q_r}\,y_{k_{m_{1}}}\cdots y_{k_{m_{j-s}}}\,y_{d_1}\cdots y_{d_{s-r}}\\ 
&&=P(k_1,\ldots,k_j, d_1,\ldots,d_{s-r})^{sb}\,(x_{k_{m_{1}}}\cdots x_{k_{m_{j-s}}}\, x_{n_{1}}\cdots x_{n_{s-r}})^{a-(j-s)b}\,w^{j-s}\\ \nonumber
&&\cdot\; y_{q_1}\cdots y_{q_r}\,y_{d_1}\cdots y_{d_{s-r}}\\ 
&&-(x_{n_1}\cdots x_{n_{s-r}})^a\,(x_{d_1}\cdots x_{d_{s-r}})^{(j-s)b}\,Q\, y_{d_1}\cdots y_{d_{s-r}}.
\end{eqnarray}

Now, expression numbered (14) is same as expression numbered (13),
but expression (15) still has way to go. In the subsequent steps we resort to Koszul generators as tags, namely, firstly,

\begin{eqnarray}\nonumber
\lefteqn{
(x_{n_2}\cdots x_{n_{s-r}})^a\,\underbrace{x_{d_0}^{a+(j-s)b}}_{=1}(x_{d_1}\cdots x_{d_{s-r}})^{(j-s)b}\,Q\, y_{d_2}\cdots y_{d_{s-r}}\,\underbrace{y_{n_0}}_{=1}\, K_{n_1,d_1} }\kern1cm \\ 
&&=(x_{n_1}\cdots x_{n_{s-r}})^a\,(x_{d_1}\cdots x_{d_{s-r}})^{(j-s)b}\,Q\, y_{d_1}\cdots y_{d_{s-r}}\\ \nonumber
&&- (x_{n_2}\cdots x_{n_{s-r}})^a\,x_{d_1}^{a+(j-s)b}(x_{d_2}\cdots x_{d_{s-r}})^{(j-s)b}\,Q\, y_{d_2}\cdots y_{d_{s-r}}\,y_{n_1}
\end{eqnarray}

The procedure establishes an inductive argument by which one monomial term is canceled against a next term in an expression involving a further down Koszul form.
To obtain the final combination in terms of earlier Sylvester forms and Koszul forms, one resorts to a summation of expressions of the same type where the first summand is the expression in the last line of the last display and the last summand recovers (11).
This explains the final form of the required expression as stated.
\qed

\smallskip

\subsubsection{Degree jump}
Now $j'=j+1$.

\begin{itemize}
\item[\textbullet] $x_{s}^{a-b},~s=j+2,\ldots,n.$

\begin{eqnarray}\nonumber
x_{s}^{a-b}H_{j+1}^{1,\ldots,j,j+1}&=&(x_{1}\cdots x_{j}x_{j+1})^{a-(j+1)b}w^{j}L_{s}\\ \nonumber
& +& x_{1}^{a-jb}(x_{j+2}\cdots x_{s-1}\,\widehat{x_{s}}\, x_{s+1}\cdots x_{n})^{b}y_{s}H_{j}^{2,\ldots,j+1}\\ \nonumber
&+& (x_{j+2}\cdots x_{s-1}\,\widehat{x_{s}}\, x_{s+1}\cdots x_{n})^{(j+1)b}x_{s}^{jb}y_{2}\cdots y_{j+1}K_{1,s}.
\end{eqnarray}

\item[\textbullet] $x_{s}^{jb},~s=1,\ldots,j+1.$

$$x_{s}^{jb}H_{j+1}^{1,\ldots,j,j+1}=(x_{1}\cdots x_{s-1}\,\widehat{x_{s}}\, x_{s+1}\cdots x_{j+1} )^{a-(j+1)b}w^{j}L_{s} - (x_{j+2}\cdots x_{n} )^{b}y_{s}H_{j}^{1,\ldots,\widehat{s},\ldots,j+1}. $$

\item[\textbullet] $(x_{r_{1}}\cdots x_{r_{s}})^{(j+1-s)b}$, where $s\in\{1,\ldots,j\}$ and $\{r_{1}, \ldots, r_{s}\}$ is an ordered subset of $\{1,\ldots,j,j+1\}.$

\begin{eqnarray}\nonumber
(x_{r_{1}}\cdots x_{r_{s}})^{(j+1-s)b}H_{j+1}^{1,\ldots,j,j+1}&=&P_{j+1}(r_1,\ldots,r_s)^{a-(j+1)b}w^{j+1-s}H_{s}^{r_{1},\ldots,r_{s}}\\ \nonumber
&+& P(1,\ldots,j+1)^{sb}y_{r_{1}}\cdots y_{r_{s}}H_{j+1-s}^{1,\ldots,\widehat{r_{1}},\ldots,\widehat{r_{s}},\ldots,j,j+1},
\end{eqnarray}
where the lower index $j+1$ in the first $P$ indicates that the product of the variables is over the complement of $\{r_{1}, \ldots, r_{s}\}$ in $\{1,\ldots,j,j+1\}$ (and not in the entire $\{1,\ldots,n\}$.)

\item[\textbullet] $(x_{q_{1}}\cdots x_{q_{r}})^{(j+1-s)b}(x_{d_{1}}\cdots x_{d_{s-r}})^{a-sb}$.

\end{itemize}

One applies the hypotheses and the conclusion of Lemma~\ref{last_case} with the following changes in the numerology:

\smallskip

\noindent $\{k_1,\ldots.k_j\} \rightsquigarrow \{1,\ldots,j+1\}$


\noindent $\{k_{m_{1}},\ldots,k_{m_{j-s}}\}\rightsquigarrow\{m_{1},\ldots,m_{j+1-s}\}$

\noindent $j \rightsquigarrow j+1$, in all appearances of $j$ in a subscript or exponent.

\medskip

Finally, we stress the calculation when the degree jumps to the highest degree in the sequence of Sylvester forms.
Once more we only display the case where $a\leq nb, \,\, (p-1)b<a\leq pb$, since when $a>nb$ the result is embedded in the general discussion of this case.

\begin{itemize}
\item[\textbullet] $x_{r}^{a-b},~r=1,\ldots,p.$\\
$$x_{r}^{a-b}H_{p}^{1,\ldots,p}=w^{p-1}L_{r}+(x_{p+1}\cdots x_{n})^{b}(x_{1}\cdots\widehat{x_{r}}\cdots x_{p})^{pb-a}y_{r}H_{p-1}^{1,\ldots,\widehat{r},\ldots,p}.$$
\item[\textbullet] $x_{s}^{a-b},~s=p+1,\ldots,n.$\\
$$x_{s}^{a-b}H_{p}^{1,\ldots,p}=w^{p-1}L_{s}+(x_{1}\cdots x_{p-1})^{pb-a}(x_{p}x_{p+1}\cdots\widehat{x_{s}}\cdots x_{n} )^{b}y_{s}H_{p-1}^{1,\ldots,p-1}$$$$\hspace{2cm}+x_{s}^{(p-1)b}(x_{1}\cdots x_{p})^{pb-a}(x_{p+1}\cdots\widehat{x_{s}}\cdots x_{n} )^{pb}y_{1}\cdots y_{p-1}K_{p,s}. $$
\item[\textbullet] $(x_{1}\cdots \widehat{x_{r}}\cdots x_{p})^{a-(p-1)b},~r=1,\ldots,p.$\\
$$(x_{1}\cdots \widehat{x_{r}}\cdots x_{p})^{a-(p-1)b}H_{p}^{1,\ldots,p}=wH_{p-1}^{1,\ldots,\widehat{r},\ldots,p}+x_{r}^{pb-a}(x_{p+1}\cdots x_{n})^{(p-1)b}y_{1}\cdots\widehat{y_{r}}\cdots y_{p}L_{r} $$


\item[\textbullet] $(x_{q_{1}}\cdots x_{q_{s}} )^{a-sb}, \,\{q_{1}<\cdots<q_{s}\}\subset\{1,\ldots,p\}$.
\begin{eqnarray}\nonumber
\lefteqn{(x_{q_{1}}\cdots x_{q_{s}} )^{a-sb}H_{p}^{1,\ldots,p}=w^{p-s}H_{s}^{q_{1},\ldots,q_{s}}}
\\\nonumber
&+&(x_{1}\cdots \widehat{x_{q_{1}}}\cdots \widehat{x_{q_{s}}}\cdots x_{p})^{pb-a}(x_{p+1}\cdots x_{n})^{sb}  y_{q_{1}}\cdots y_{q_{s}}H_{p-s}^{1,\ldots,\widehat{q_{1}},\ldots,\widehat{q_{s}}, \ldots, p}
\end{eqnarray}

\item[\textbullet] $(x_{q_{1}}\cdots x_{q_{r}}x_{d_{1}}\cdots x_{d_{s-r}})^{a-sb}$, where $s\in\{2,\ldots,p-1\},  r\in\{0,\ldots,s-1\},\{q_{1}<\ldots <q_{r}\}\subset\{1,\ldots,p\},  p<d_{1}<\cdots<d_{s-r}$.

This case follows the pattern of Lemma~\ref{last_case}, with the following changes:

\smallskip

\noindent $\{k_1,\ldots, k_j\} \rightsquigarrow \{1,\ldots,p\}$

\noindent $\{m_{1},\ldots,m_{j-s}\}\rightsquigarrow\{m_{1},\ldots,m_{p-s}\}$

\noindent $j \rightsquigarrow p$, in all appearances of $j$ in a subscript or exponent.

\smallskip

However, there are some changes in the coefficients of the final expression.
We write this expression for the sake of completeness:

\begin{eqnarray}\nonumber
\lefteqn{(x_{q_{1}}\cdots x_{q_{r}}x_{d_{1}}\cdots x_{d_{s-r}})^{a-sb}H_{p}^{1,\ldots,p}=w^{p-s}H_{s}^{q_{1},\ldots,q_{r},d_{1},\ldots,d_{s-r}} }\kern0.5cm\\\nonumber
&& +P(q_1,\ldots,q_r,m_1,\ldots,m_{p-s}, d_1,\ldots, d_{s-r})^{sb}(x_{m_{1}}\cdots x_{m_{p-s}} )^{pb-a}\\\nonumber
&&\cdot
y_{q_{1}}\cdots y_{q_{r}}y_{d_{1}}\cdots y_{d_{s-r}}H_{p-s}^{m_{1},\ldots,m_{p-s}} 
+ \sum_{c=1}^{s-r}(x_{p+1}\cdots \widehat{x_{d_{1}}}\cdots \widehat{x_{d_{s-r}}} \cdots x_{n}  )^{pb}  \\\nonumber
&&\cdot (x_{d_{c-1}}\cdots x_{d_{s-r}} )^{a} (x_{d_{1}}\cdots x_{d_{s-r}})^{(p-s)b}(x_{m_{1}}\cdots x_{m_{p-s}})^{pb-a}\\\nonumber
&&\cdot  (x_{n_{c+1}}\cdots x_{n_{s-r}})^{pb}(x_{n_{1}}\cdots x_{n_{c}})^{pb-a}(x_{q_{1}}\cdots x_{q_{r}})^{(p-s)b}
\\\nonumber
&&\cdot y_{q_{1}}\cdots y_{q_{r}}y_{m_{1}}\cdots y_{m_{p-s}}y_{n_{1}}\cdots y_{n_{c-1}}y_{d_{c+1}} \cdots y_{d_{s-r}}K_{n_{c},d_{c}},
\end{eqnarray}
with the convention that $x_{d_0}=y_{n_0}=x_{n_{s-r+1}}=x_{d_{s-r+1}}=1$.

\item[\textbullet] $(x_{1}\cdots x_{j}\, x_{r_{j+1}}\cdots x_{r_{p-1}})^{a-(p-1)b}$, with $j\in\{0,\ldots,p-2\},~k\in \{j+1,\ldots,p-1\}$ and $r_{k}\in \{p+1,\ldots,n\}.$
\end{itemize}

This case still follows the general shape afforded by the result of Lemma~\ref{last_case}, except that, besides changes in the $\mathbf{x}$-coefficients, the Sylvester tag of one of the terms degenerates into a syzygy generator $L_{j+1}$.
\begin{eqnarray}\nonumber
\lefteqn{\kern 0.3cm (x_{1}\cdots x_{j}x_{r_{j+1}}\cdots x_{r_{p-1}})^{a-(p-1)b}H_{p}^{1,\ldots,p}=  wH_{p-1}^{1,\ldots,j,r_{j+1},\ldots, r_{p-1}}} \kern0.1cm\\\nonumber
&& +x_{j+1}^{pb-a}(x_{j+2}\cdots \widehat{x_{r_{j+1}}}\cdots \widehat{x_{r_{p-1}}}\cdots x_{n} )^{(p-1)b}\, y_{1}\cdots y_{j}y_{r_{j+1}}\cdots y_{r_{p-1}}L_{j+1}\\\nonumber
&&+\sum_{c=2}^{p-1-j}(x_{1}\cdots x_{j})^{b}(x_{j+1}\cdots x_{j+c})^{pb-a}(x_{j+c+1}\cdots x_{p})^{pb}(x_{r_{j+1}}\cdots x_{r_{j+c-2}})^{a+b}(x_{r_{j+c-1}}\cdots x_{r_{p-1}} )^{b} \\\nonumber
&& \cdot (x_{p+1}\cdots \widehat{x_{r_{j+1}}}\cdots \widehat{x_{r_{p-1}}} \cdots x_{n})^{pb}y_{1}\cdots y_{j+c-1}y_{r_{j+c}}\cdots y_{r_{p-1}}K_{j+c,r_{j+c-1}}\\\nonumber
&& +(x_{1}\cdots x_{j})^{b}(x_{j+1}\cdots x_{p})^{pb-a}(x_{r_{j+1}}\cdots x_{r_{p-2}})^{a+b}x_{r_{p-1}}^{b} (x_{p+1}\cdots \widehat{x_{r_{j+1}}}\cdots \widehat{x_{r_{p-1}}} \cdots x_{n})^{pb} \\\nonumber
&& \cdot y_{1}\cdots y_{p-1}K_{p,r_{p-1}}.
\end{eqnarray}

This concludes the proof of the proposition.
\qed


                       \end{document}